\documentclass[11pt, a4paper]{article}

\usepackage{graphicx}
\usepackage{color}
\usepackage{amsmath}
\usepackage{amssymb}
\usepackage{amscd}
\usepackage{bbm}
\usepackage{tikz}
\usetikzlibrary{patterns}
\usepackage{hyperref}
\hypersetup{
    bookmarks=true,         
    colorlinks=true,       
    linkcolor=red,          
    citecolor=green,        
    filecolor=magenta,      
    urlcolor=cyan           
}

\usepackage[utf8]{inputenc}
\usepackage{amsthm}
\usepackage{bbm} 
\usepackage[linesnumbered,lined,boxed,commentsnumbered]{algorithm2e}

\usepackage{marginnote}

\newtheorem{Assumption}{Assumption}
\newtheorem{Definition}{Definition}

\newtheorem{Remark}{Remark}
\newtheorem{Example}{Example}

\newcommand{\inr}[1]{\bigl< #1 \bigr>}

\newcommand{\norm}[1]{\left\|#1\right\|}%
\newcommand{\an}{{\mbox{ and }}}

\newcommand{\pa}[1]{\left(#1\right)}
\newcommand{\cro}[1]{\left[#1\right]}
\newcommand{\set}[1]{\left\{#1\right\}}

\newcommand\eps{\epsilon}

\DeclareMathOperator*{\argmin}{argmin}

\newcommand{\ERM}[1]{\widehat{f}_{#1}}

\def\ds1{\textrm{1\kern-0.25emI}} 

\newcommand{\1}{\ensuremath{\mathbbm{1}}}

\newcommand \E{\mathbb{E}}
\newcommand \R{\mathbb{R}}

\newcommand \cC{{\cal C}}
\newcommand \cD{{\cal D}}

\newcommand \cF{{\cal F}}
\newcommand \cG{{\cal G}}

\newcommand \cI{{\cal I}}
\newcommand \cJ{{\cal J}}
\newcommand \cK{{\cal K}}
\newcommand \cL{{\cal L}}

\newcommand \cN{{\cal N}}
\newcommand \cO{{\cal O}}

\newcommand \cR{{\cal R}}
\newcommand \cS{{\cal S}}

\newcommand \cX{{\cal X}}
\newcommand \cY{{\cal Y}}

\newcommand \bN{{\mathbb N}}

\newcommand \bP{{\mathbb P}}

\newcommand{\bayes}{f^*}
\newcommand{\MOM}[2]{\text{MOM}_{#1}\big(#2\big)}

\newcommand{\param}[1]{\theta_{#1}}

\usepackage{amsmath, amsfonts, amssymb, amsthm, stmaryrd} 
\usepackage{enumitem}
\usepackage[lofdepth,lotdepth]{subfig}
\usetikzlibrary{arrows}

\renewcommand{\P}{\mathbb{P}}

\newtheorem{Prop}{Proposition} 
\newtheorem{Th}{Theorem}

\newcommand{\llb}{\left\llbracket}
\newcommand{\rrb}{\right\rrbracket}

\begin{document}

\title{Robust classification via MOM minimization}
\author{G. Lecué and M. Lerasle and T. Mathieu}
\maketitle
\begin{abstract}
We present an extension of Vapnik's classical empirical risk minimizer (ERM) where the empirical risk is replaced by a median-of-means (MOM) estimator, the new estimators are called MOM minimizers. While ERM is sensitive to corruption of the dataset for many classical loss functions used in classification, we show that MOM minimizers behave well in theory, in the sense that it achieves Vapnik's (slow) rates of convergence under weak assumptions: data are only required to have a finite second moment and some outliers may also have corrupted the dataset. 

We propose an algorithm inspired by MOM minimizers. 
These algorithms can be analyzed using arguments quite similar to those used for Stochastic Block Gradient descent. 
As a proof of concept, we show how to modify a proof of consistency for a descent algorithm to prove consistency of its MOM version. 
As MOM algorithms perform a smart subsampling, our procedure can also help to reduce substantially time computations and memory ressources when applied to non linear algorithms.
These empirical performances are illustrated on both simulated and real datasets.
\end{abstract}

\section{Introduction}

Machine learning has experienced unprecedented growth in recent years with a major societal impact exceeding that it has had in science for several decades.
The article is set within the framework of robust learning theory, in the sense that only  moment hypotheses are guaranteed on the data and the dataset may contain outliers.
Robust learning theory has received particular attention in recent years, the aim sometimes being to construct reliable procedures under weaker assumptions than strong stochastic ones (such as the Gaussian character of noise in regression) and sometimes to detect outliers.
This interest can be appreciated, for example, by the challenges recently posted on ``kaggle", the most popular data science competition platform. The 1.5 million dollars problem ``Passenger Screening Algorithm Challenge'' involves the discovery of terrorist activity from 3D images. The ``NIPS 2017: Defense Against Adversarial Attack'' consists in building algorithms robust to adversarial data.

\begin{figure}[h]
\begin{center}
\includegraphics[scale=0.55]{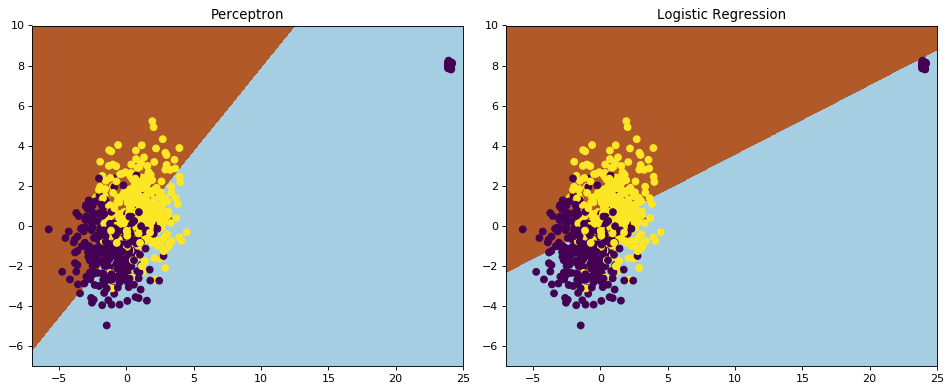}
\caption{Scatter plot of the toy dataset, the color of the points gives their class. The background color gives the linear separation provided by the perceptron (left) and the logistic regression (right) trained on this corrupted dataset.\label{fig:intro}}
\end{center}
\end{figure}

Our approach extends Vapnik's \emph{empirical risk minimization} (ERM) (cf., \cite{MR1641250}). The ERM is defined as a minimizer of the empirical risk $f\to (1/N)\sum_{i=1}^NI_{\{Y_i\ne f(X_i)\}}$ over a class $\cF$ of functions from $\cX$ to $\{\pm1\}$:
\begin{equation}\label{eq:DefERM}
\ERM{\text{ERM}} \in \argmin_{f\in \cF}\frac1N\sum_{i=1}^NI{\{Y_i\ne f(X_i)\}} 
\end{equation}where $I{\{Y_i\ne f(X_i)\}} = 1$ if $Y_i\ne f(X_i)$ and $0$ otherwise.
This estimator is fairly insensitive to both heavy-tailed data and to the presence of a small proportion of outliers in the dataset (see Theorem~\ref{Lem:ERMOutliers} below). The reason is that the $0-1$ loss $\ell_f(x,y)=I_{\{y\ne f(x)\}}$ is bounded, which grants concentration no matter the distribution of $X$ and a small number of data cannot really impact the empirical mean performance. 
However, in practice, we find that algorithms supposed to approach ERM (usually following a convex relaxation principle) have poor performance in the presence of a few outliers.
Figure~\ref{fig:intro} illustrates this problem on a toy example on which most of the data would be well separated by a linear classifier like Perceptron \cite{MR0122606} or logistic classifier, but adding some anomalies completely flaws these algorithms.

The example in Figure~\ref{fig:intro} is representative of a general problem. 
Algorithms approaching ERM do not directly solve the minimization problem \eqref{eq:DefERM}, which is NP-hard even for simple function classes like half-spaces indicators \cite{MR2529778,MR3029261}. 
Instead, the algorithm seeks to approach a convex relaxation of the problem, in which the considered loss (e.g. hinge or logistic loss) is not necessarily bounded. 
The advantage of this relaxation is that the function to be minimized is convex and that the minimum can easily be approached using for example gradient descent algorithms.
On the other hand, this function becomes unbounded and can be diverted by a single outlier.
Minimizing this corrupted criterion then leads to a bad classifier choice.

In learning theory, most alternatives to ERM manage the problem of outliers and heavy tail distributions for outputs only. These solutions are based on the pioneering work of John Tukey \cite{tukey1960survey,tukey1962future}, Peter Huber \cite{MR0161415,huber1967behavior} and Frank Hampel \cite{hampel1971general,hampel1974influence}, replacing the square loss by a robust alternative like Huber loss or Tuckey's biweight loss.
These methods do not allow to treat the case where the inputs are with heavy tails or corrupted, which is a classical problem in robust statistics also known as the ``\textit{leverage point problem}'', see \cite{robuststat}. In this article, we address this question by considering an alternative to M-estimators, called median-of-means. 
Several estimators based on this principle have recently been proposed in the literature \cite{MR3378468,MOM1,LugosiMendelson2017,LugosiMendelson2017-2,LugosiMendelson2016,shahar_mom_1,MOM2}. 
To our knowledge, these articles use the small ball hypothesis \cite{MR3431642,Shahar-COLT} to treat problems of least square regression or Lipschitzian loss regression.
This assumption is restrictive in some classic functional regression frameworks \cite{MR3757527, Wellner:2017} or for problems such as the construction of recommendation system where inputs are sampled in the canonical basis and therefore do not satisfy a small ball condition. In this article, we propose a natural estimator based on the MOM principle, which we call MOM minimizer, and which we study without the small ball hypothesis. Instead we assume an a priori bound on the $L^2$-norm of learning functions. We can identify mainly two streams of hypothesis in Learning theory: 1) boundedness with respect to some norm of the class $F$ of functions and the output $Y$, the typical example is the boundedness in $L_\infty$ assumption or 2) norm equivalence assumption over the class $F$ (or, more precisely, on the shifted class $F-f^* = \{f-f^*:f\in F\}$ where $f^*$ is the oracle in $F$) and $Y$, the typical example being the subgaussian assumption, i.e. $\norm{f-f^*}_{\psi_2}\leq L \norm{f-f^*}_{L_2}, \forall f\in F$.  The small ball assumption is a norm equivalence assumption between the $L_1$ and $L_2$ norms and is concerned with the second type of assumptions. Our approach deals with the first type of assumption since we only assume boundedness in $L_2$ which can be seen as an significant improvement upon a $L_\infty$ boundedness assumption even though the rates of convergence obtained for our estimators turn out to be minimax \cite{probapatern} in the absence of a margin condition \cite{MR1765618} even under a $L_\infty$ assumption.

The estimation of the expectation of a univariate variable by median-of-means (MOM) \cite{MR1688610, MR855970, MR702836} is done as follows: given a partition of the dataset into blocks of the same size, the empirical mean is constructed on each block and the MOM estimator is obtained by taking the median of these empirical means (see Section~\ref{sec:CorrDataSets} for details). These estimators are naturally resistant to the presence of a few outliers in the dataset: if the number of these outliers does not exceed half the number of blocks, more than half of these blocks are made of ``clean" data and the median is a reliable estimator. 

In this article, we also study a family of algorithms approaching the MOM minimizer.
The general idea we wish to promote is to use the MOM principle within algorithms originally intended for the evaluation of ERM associated to convex loss functions.
In Section~\ref{sec:Comp}, we present the modification of gradient descent algorithms following this philosophy. 
The general principle of this iterative algorithm is as follows: 
at iteration $t$, a dataset equipartition $B_1,\ldots,B_K$ is selected uniformly at random and the most central block $B_{\text{med}}$ is determined according to the following formula
\begin{equation}\label{eq:MOM_estimate}
 \sum_{i\in B_{\text{med}}}\ell_{f_t}(X_i,Y_i)=\text{median}\pa{\sum_{i\in B_{k}}\ell_{f_t}(X_i,Y_i):\quad k=1,\ldots,K}\enspace.
\end{equation}
The next iteration $f_{t+1}$ is then produced by taking from $f_t$ a step down in the direction opposite to the gradient of  $f \to \sum_{i\in B_{\text{med}}}\ell_{f}(X_i,Y_i)$ at $f_t$, cf. Algorithm~\ref{fig:mom_algo}. 
The underlying heuristic is that the data in the selected block $B_{\text{med}}$  are safe for estimating the mean at the point, in the sense that empirical risk $|B_{{\rm med}}|^{-1}\sum_{i\in B_{\text{med}}}\ell_{f_t}(X_i,Y_i)$ is a subgaussian estimator of $\E \ell_{f}(X_i,Y_i)$, cf. \cite{MR3576558} and that data indexed by $B_{\text{med}}$ should not be outliers. The consistency of the algorithm is studied in Section~\ref{sec:conv_MOM_algo}. 
Theorem~\ref{thm:convergence} provides an example showing how to modify the convergence proof of SGD algorithms.

In Section~\ref{simu}, the practical performances of this algorithm are illustrated on several simulations, involving in particular different loss functions. 
These simulations illustrate not surprisingly the gain of robustness that there is to use these algorithms in their MOM version rather than in their traditional version, as can for example be appreciated on the toy-example of Figure~\ref{fig:intro}, cf. Figure~\ref{fig:classif_comparaison}. 
MOM estimators are compared to different learning algorithms on real datasets that can be modeled by heavy tailed data, obtaining in each case performances comparable to the best of these benchmarks.

Another advantage of our procedure is that it works on blocks of data. 
This can improve speed of execution and reduce memory requirements, which can be decisive on massive datasets and/or when one wishes to use non-linear algorithms as in Section~\ref{complexity}. 
This principle of dividing the dataset to calculate estimators more quickly and then aggregating them is a powerful tool in statistics and machine learning \cite{Jor:2013}.
Among others, one can mention bagging methods \cite{Bre:1996a} or subagging ---a variant of bagging where the bootstrap is replaced by subsampling--- \cite{Buh_Yu:2002}. 
These methods are considered difficult to study theoretically in general and their analysis is often limited to obtaining  asymptotic guarantees.
By contrast, the theoretical tools for non-asymptotic risk analysis of MOM minimizer have already essentially been developed and the consistency of the algorithm can also be derived from standard arguments.
Finally, subsampling by the central block $B_{{\rm med}}$ ensures robustness properties that cannot be guaranteed by traditional alternatives.  

The algorithm provides an empirical notion of data depth: 
data providing good risk estimates of $f\to \E\ell_f(X,Y)$ are likely to be selected many times in the central block $B_{\text{med}}$ along the descent, while outliers will be systematically ignored. 
This notion of depth, based on the risk function, is very natural for prediction problems.
It is complemented by an outliers detection procedure: data that are selected a number of times below a predetermined threshold are classified outliers. 
This procedure is evaluated on our toy example of Figure~\ref{fig:intro}. 
Data represented by the dots in the top right corner all end with a null score (see Figure~\ref{fig:outlier1} below).
The procedure is then tested on a real dataset on which the conclusions are more interesting.
On this experiment, according to the theoretical bounds in Theorem~\ref{thm:MOM}, MOM minimizer's prediction qualities are all the better when the number of blocks is small, cf. Figure~\ref{fig:boxplot_pulsar}.
On the other hand, outlier detection is best when the number of blocks is large, cf. Figure~\ref{fig:outlier_pulsar}.
Outlier resistance and anomaly detection tasks can therefore both be handled using the MOM principle, but the main hyper-parameter $K$ -- the number of blocks of data --  for setting this method must be chosen carefully according to the objective.
A number of blocks as small as possible (about twice the number of outliers) will give the best predictions, while large values of this number of blocks will accelerate the detection of anomalies.
Note that it is essential for outliers detection to use different (for instance, random) partitions at each step of the descent to avoid giving the same score to an outlier and to all the data in the same block containing it.

Detecting outliers is usually performed in machine learning via some unsupervised preprocessing algorithm that detects outliers outside a bulk of data, see for example \cite{robust_LDA,robust_LDA_2,robust_LDA_3,robust_LDA_4} or other algorithms like DBSCAN \cite{birant2007st} or isolation forest \cite{liu2008isolation}. 
These algorithms assume elliptical symmetry of the data, a solution for skewed data can also be found in \cite{Hubert_skewed}. 
These unsupervised preprocessing removes outliers in advance, i.e. before starting any learning task. 
As expected, these strategies work well in the toy example from Figure~\ref{fig:intro}. 
There are several cases where it will fail though. 
First, as explained in \cite{robuststat}, this strategy classifies data independently of the risk, it is likely to remove from the dataset outlier coming from heavy-tailed distribution, yielding biased estimators. 
Moreover, a small group of misclassified data inside a bulk won't be detected. 
Our notion of depth, based on the risk, seems more adapted to the learning task than any preprocessing procedure blind to the risk.

The paper is organized as follows. Section~\ref{sec:Setting} presents the classification problem, the ERM and its MOM versions and gathers the assumptions granted for the main results. Section~\ref{sec:MainThRes} presents theoretical risk bounds for the ERM and MOM minimizers on corrupted datasets. Section~\ref{sec:Comp} deals with theoretical results on the algorithm computing MOM minimizers. We present the algorithm, study its convergence and provide theoretical bounds on its complexity. Section~\ref{simu} shows empirical performance of our estimators in both simulated and real datasets. Proofs of the main results are postponed to Section~\ref{sec:Proof} where we also added heuristics on the practical choice of the hyper-parameters.

\section{Setting}\label{sec:Setting}
\subsection{Empirical risk minimization for binary classification}\label{sec:Classif}
Consider the supervised binary classification problem, 
where one observes a sample $(X_1,Y_1)$, $\ldots$, $(X_N,Y_N)$ taking values in $\cX \times \cY$. The set $\cX$ is a measurable space and $\mathcal{Y} = \{-1,1\}$. 
The goal is to build a classifier ---that is, a measurable map $f:\cX \to \cY$--- 
such that, for any new observation $(X,Y)$, $f(X)$ is a good prediction for $Y$. 
For any classifier $f$, let
\[
\ell^{0-1}_f(x,y)=I\{y\ne f(x)\},\qquad R^{0-1}(f) = P\ell^{0-1}_f=\bP_{(X,Y)\sim P} \bigl( Y\neq f(X) \bigr)\enspace.
\]
The $0-1$ risk $R^{0-1}(\cdot)$ is a standard measure of the quality of a classifier.
Following Vapnik \cite{MR1719582}, a popular way to build estimators is to replace the unknown measure $P$ in the definition of the risk by the empirical measure $P_N$ defined for any real valued function $g$ by $P_Ng=N^{-1}\sum_{i=1}^Ng(X_i,Y_i)$ and minimize 
the empirical risk. The \emph{empirical risk minimizer} for the $0-1$ loss on a class $\cF$ of classifiers is $\ERM{\text{ERM}}^{0-1}\in\argmin_{f\in \cF}\{P_N\ell^{0-1}_f\}$.

The main issue with $\ERM{\text{ERM}}^{0-1}$ is that it cannot be computed efficiently in general. One source of computational complexity is that both $\cF$ and the $0-1$ loss function are non-convex. This is why various convex relaxations of the $0-1$ loss have been introduced since the very beginning of machine learning. These proceed in two steps. First, $\cF$ should be replaced by a convex set $F$ of functions taking values in $\R$. Then one builds an alternative loss function $\ell$ for $\ell^{0-1}$ defined for all $f\in F$. The new function $\ell$ should be convex and put less weight on those $f\in F$ such that $f(X_i)Y_i>0$. Classical examples include the\emph{ hinge loss} $\ell^{\text{hinge}}_f(x,y)=(1-yf(x))_+$, or the \emph{logistic loss} $\ell^{\text{logistic}}_f(x,y)=\log(1+e^{-yf(x)})$. A couple $(F,\ell)$ such that $F$ is a convex set of real valued functions and $\ell$ is a convex function (i.e. for all $y\in\{-1, 1\}$ and $x\in\cX$, $f\in F\to \ell_f(x, y)$ is convex) such that $\ell_f(x,y)<\ell_f(x,-y)$ whenever $yf(x)>0$ will be called a convex relaxation of $(\cF,\ell^{0-1})$. Given a convex relaxation $(F,\ell)$ of $(\cF,\ell^{0-1})$, one can define the associated empirical risk minimizer (ERM) by
\begin{equation}\label{def:ERM}
 \ERM{\text{ERM}}\in\argmin_{f\in F} P_N\ell_f\enspace.
\end{equation} Note that $\ERM{\text{ERM}}$ does not build a classifier. To deduce, a classification rule from $\ERM{\text{ERM}}$ one can simply consider its sign function defined for all $x\in \cX$ by ${\rm sign}(\ERM{\text{ERM}}(x)) = 2(I\{\ERM{\text{ERM}}(x)\ge 0\}-1/2)$. The procedure $ \ERM{\text{ERM}}$ is solution of a convex optimization problem that can therefore be approximated using a descent algorithm. We refer for example to \cite{BubeckOptim} for a recent overview of this topic and Section~\ref{sec:Comp} for more examples.

\subsection{Corrupted datasets}\label{sec:CorrDataSets}

In this paper, we consider a framework where the dataset may have been corrupted by \emph{outliers} (or anomalies). There are several definitions of outliers in the literature, here, we assume that the dataset is divided into two parts. The first part is the set of inliers, indexed by $\cI$, data $(X_i,Y_i)_{i\in \cI}$ are hereafter always assumed to be independent and identically distributed (i.i.d.) with common distribution $P$. The second one is the set of outliers, indexed by $\cO\subset[N]$ which has cardinality $|\cO|$. Nothing is assumed on these data which may not be independent, have distributions $P_i$ totally different from $P$, with no moment at all, etc... In particular, this framework is sufficiently general to cover the case where outliers are i.i.d. with distribution $Q\neq P$ as in the $\eps$-contamination model  \cite{robuststat,chao_1,chao_2,montanari_donoho}.

Our first result shows that the rate of convergence of $\ERM{\text{ERM}}^{0-1}$ is not affected by this corruption as long as $|\cO|$ does not exceed $N\times(\text{rate of convergence})$ see Theorem~\ref{Lem:ERMOutliers} and the remark afterward. However, it is easy to remark that, when the number $N$ of data is finite as it is always the case in practice, even one aggressive outliers may yield disastrous breakdown of the empirical mean's statistical performance. Consequently, even if $\ERM{\text{ERM}}^{0-1}$ behaves correctly, its proxy $\ERM{\text{ERM}}$ defined in \eqref{def:ERM} for a convex relaxation $(F, \ell)$ can have disastrous statistical performances, particularly when $F$ and $\ell$ are unbounded, cf. Figure~\ref{fig:classif_comparaison} for an illustration. 

To bypass this problem, we consider in this paper an alternative to the empirical mean called \emph{median-of-means} \cite{MR1688610, MR855970, MR702836}. Let $K\leq N$ denote an integer and let $B_1,\ldots,B_K$ denote a partition of $\{1,\ldots,N\}$ into bins $B_k$ of equal size $|B_k|=N/K$. If $K$ doesn't divide $N$, one can always drop a few data. For any function $f:\cX\times\cY\to\R$ and any non-empty subset $B\subset \{1,\ldots,N\}$, define the empirical mean on $B$ by $P_Bf=|B|^{-1}\sum_{i\in B}f(X_i,Y_i)$. The median-of-means (MOM) estimator of $Pf$ is defined as the empirical median of the empirical means on the blocks $B_k$
\[
\MOM{K}{f}=\text{median}\set{P_{B_k}f:\;k=1,\ldots,K}\enspace.
\]
As the classical Huber's estimator \cite{MR0161415}, MOM estimators interpolate between the unbiased but non robust empirical mean and the robust but biased median. 
%
%
%
%
In particular, when applied to loss functions, these new estimators of the risk $P\ell_f, f\in F$ suggest to define the following alternative to Vapnik's ERM, called MOM minimizers
\begin{equation}\label{eq:defMOMK}
\ERM{\text{MOM},K}\in\argmin_{f\in F} \MOM{K}{\ell_f}\enspace. 
\end{equation}
From a theoretical point of view, we will prove that, when the number $|\cO|$ of outliers is smaller than $N\times(\text{rate of convergence})$, $\ERM{\text{MOM},K}$ performs well under $2$-moments assumptions on $F$ and $\ell$. To illustrate our main assumptions and theoretical results, we will regularly use the following classical example.

\begin{Example}[linear classification.]\label{Ex:linClass}
Let $\cX=\R^p$ and let $\norm{\cdot}_2$ denote the classical Euclidean norm on $\R^p$. Let $F$ denote a set of linear functions
\[
F=\{f_t:x\mapsto \inr{x, t}: \norm{t}_{2}\le \Gamma\}\enspace.
\]
Let $\ell$ denote either the hinge loss or the logistic loss defined respectively for any $(x,y)\in\cX\times\cY$ and $f\in F$ by
\[
\ell^{\text{hinge}}_f(x,y)=(1-yf(x))_+,\qquad \ell^{\text{logistic}}_f(x,y)=\log(1+e^{-yf(x)})\enspace.
\]

\end{Example}


\subsection{Main assumptions}
As already mentioned, data are divided into two groups, a subset $\{(X_i,Y_i) : i\in\cO\}$ made of outliers (on which we will make no assumption) and the remaining data $\{(X_i,Y_i): i\in\cI\}$ contains all data that bring information on the target/oracle 
\[
\bayes\in\argmin_{f\in F}P\ell_f\enspace.
\]
Data indexed by $\cI$ are therefore called \emph{inliers} or \emph{informative data}. To keep the presentation as simple as possible, inliers are assumed to be i.i.d. distributed according to $P$ although this assumption could be relaxed as in \cite{MOM1,MOM2}. Finally, note that the $\cO\cup\cI = \{1, \ldots, N\}$ partition of the dataset is of course unknown from the statistician.

Let us now turn to the set of assumptions we will use to study MOM minimizers procedures. For any measure $Q$ and any function $f$ for which it makes sense, denote by $Qf=\int f dQ$. Denote also, for all $q\geq1$, by $L^q$ the set of real valued functions $f$ such that $\int |f|^qdP<\infty$ and, for any $f\in L^q$, by 
\[
\norm{f}_{L^q}=\pa{\int |f|^qdP}^{1/q}\enspace.
\]
Our first assumption is an $L^2$-assumption on the functions in $F$.
\begin{Assumption}\label{Ass:L2}
For all $f \in F$, we have $\norm{f}_{L^2}\le \param{2}$.
\end{Assumption}
Of course, Assumption~\ref{Ass:L2} is granted if $F$ is a set of classifiers. It also holds for the linear class of functions from  Example~\ref{Ex:linClass} as long as $P\norm{X}_2^2<\infty$ with $\param{2}=\Gamma(P\norm{X}_2^2)^{1/2}$.

The second assumption deals with the complexity of the class $F$. This complexity appears in the upper bound of the risk. It is defined using only informative data. Let 
\[
\cK=\{k\in\{1,\ldots,K\} : B_k\cap\cO=\emptyset\} \quad\an\quad \cJ=\cup_{k\in\cK}B_k\enspace.
\]
\begin{Definition}\label{Def:Rad}
Let $\cG$ denote a class of functions $f:\cX\to\R$ and let $(\epsilon_i)_{i\in \cI}$ denote i.i.d. Rademacher random variables independent from $(X_i, Y_i)_{i\in \cI}$. The \emph{Rademacher complexity} of $\cG$ is defined by
\[
\cR(\cG)=\max_{A\in\{\cI,\cJ\}}\E\left[\sup_{f\in \cG}\sum_{i\in A}\epsilon_if(X_i)\right]\enspace.
\] 
\end{Definition}
The Rademacher complexity is a standard measure of complexity in classification problems \cite{MR1984026}. It can be upper bounded, depending on the situation, by the VC dimension or the Dudley's entropy integral or the Gaussian mean width of the class $F$, see for example \cite{MR2182250,MR2829871,MR1984026,concentration,probapatern} for a presentation of these classical bounds. Our second assumption is simply that the Rademacher complexity of the class $F$ is finite.
\begin{Assumption}\label{Ass:RadComp}
The Rademacher complexity of $F$ is finite, $\cR(F)<\infty$. 
\end{Assumption}
Assumption~\ref{Ass:RadComp} holds in the linear classification example under Assumption~\ref{Ass:L2} since it follows from Cauchy-Schwarz inequality that $\cR(F)\le \param{2}\sqrt{|\cI|p}$.
Finally, our last assumption is that the loss function $\ell$ considered is Lipschitz in the following sense.
\begin{Assumption}\label{Ass:Lip}
 The loss function $\ell$ satisfies for all $(x,y)\in \cX\times\cY$ and all $f,f^\prime \in F$,
 \[
|\ell_f(x,y)-\ell_{f'}(x,y)|\le L|f(x)-f'(x)| \enspace.
 \]
\end{Assumption}
Assumption~\ref{Ass:Lip} holds for classical convex relaxation of the $0-1$ loss such as hinge loss $\ell^{\text{hinge}}$ or logistic loss $\ell^{\text{logistic}}$ as in Example~\ref{Ex:linClass}. In these examples, the constant $L$ can be chosen equal to $1$.

\section{Theoretical guarantees}\label{sec:MainThRes}
Our first result follows Vapnik's original risk bound for the ERM and shows that $\ERM{\text{ERM}}^{0-1}$ is insensitive to the presence of outliers in the dataset. Moreover, it quantifies this robustness property since Vapnik's rate of convergence is still achieved by $\ERM{\text{ERM}}^{0-1}$ when there are less than (number of observations) times (Vapnik's rate of convergence) outliers.

\begin{Th}\label{Lem:ERMOutliers}
Let $\cF$ denote a collection of classifiers. Let $\cL_\cF^{0-1}=\{\ell_f^{0-1}-\ell_{f^*}^{0-1} : f \in \cF \}$ be the family of excess loss functions indexed by $\cF$ where $f^*\in\argmin_{f\in \cF} R^{0-1}(f)$. For all $K>0$, with probability at least $1-e^{-K}$, we have
$$R^{0-1}(\ERM{\text{ERM}}^{0-1})-\inf_{f\in \cF}R^{0-1}(f)\le \frac{2\cR(\cL_\cF^{0-1})}{ N}+\sqrt{\frac{K}{2|\cI|}}+\frac{2|\cO|}{N}.$$
\end{Th}

Theorem~\ref{Lem:ERMOutliers} is proved in Section~\ref{sec:ProofLem1}. It is an adaptation of Vapnik's proof of the excess risk bounds satisfied by  $\ERM{\text{ERM}}^{0-1}$ in the presence of outliers.

\begin{Remark}
In the last result, one can easily bound the excess risk using $\cR(\mathcal{F})$ instead of $\cR(\cL_\cF^{0-1})$ since
$$\cR(\cL_\cF^{0-1}) = \max_{A\in\{\cI,\cJ\}}\E\left[\sup_{f\in \cF}\sum_{i\in A}\epsilon_i(f(X_i)-f^*(X_i))\right] = \cR(\mathcal{F}) \enspace.$$
The final bound is of similar flavor: for all $K>0$, with probability at least $1-e^{-K}$, we have
\begin{equation}\label{eq:RBERM}
 R^{0-1}(\hat f_{\text{ERM}}^{0-1})-\inf_{f\in \cF}R^{0-1}(f)\lesssim\max\left(\frac{\cR(\cF)}{N},\sqrt{\frac{K}{|\cI|}}, \frac{|\cO|}{N}\right)\enspace.
 \end{equation}
\end{Remark}

\begin{Remark}
 When $\cF$ is the  class of all linear classifiers, that is when $\cF=\{{\rm sgn}(\inr{t,\cdot}):  t\in\R^p\}$, one has $\cR(\cF)\le \sqrt{|\cI| p}$ (see Theorem~3.4 in \cite{MR2182250}). Therefore, when $|\cI|\geq N/2$, Theorem~\ref{Lem:ERMOutliers} implies that for all $1\leq K\leq p$, with probability at least $1-\exp(-K)$,
\[
R^{0-1}(\hat f_{\text{ERM}}^{0-1})-\inf_{f\in \cF}R^{0-1}(f)\lesssim\max(\sqrt{p/N}, |\cO|/N)\enspace.
\]
As a consequence, when the number of outliers is such that $|\cO|\lesssim N\times\sqrt{p/N}$, Vapnik's classical ``slow" rate  of convergence  $\sqrt{p/N}$ is still achieved by the ERM even if $|\cO|$ outliers have polluted the dataset. The interested reader can also check that ``fast rates" $p/N$ could also be achieved by the ERM in the presence of outliers if $|\cO|\lesssim p$ and when the so-called strong margin assumption holds (see, \cite{MR2182250,MR2182250}). Note also that the previous remark also holds if $F$ is a class with VC dimension $p$ beyond the case of indicators of half spaces. 
\end{Remark}

The conclusion of Theorem~\ref{Lem:ERMOutliers} can be misleading in practice. 
Indeed, theoretical performance of the ERM for the $0-1$ loss function are not downgraded by outliers, but its proxies based on convex relaxation  $(F,\ell)$  of $(\cF,\ell^{0-1})$ are. 
This can be seen on the toy example in Figure~\ref{fig:intro} and in Figure~\ref{fig:classif_comparaison} from Section~\ref{simu}. 
In this work, we propose a robust surrogate, based on MOM estimators of the risk and defined in \eqref{eq:MOM_estimate}, to the natural empirical risk estimation of the risk which works for unbounded loss functions. In the next result, we prove that the MOM minimizer $\hat f_{\text{MOM},K}$ defined as 
\begin{equation}\label{eq:MOMK_esti}
 \hat f_{\text{MOM},K}\in\argmin_{f\in F} \MOM{K}{\ell_f}
 \end{equation} satisfies an excess risk bound under weak assumptions introduced in Section~\ref{sec:Setting}. 

\begin{Th}\label{thm:MOM}
Grant Assumptions \ref{Ass:L2}, \ref{Ass:RadComp} and \ref{Ass:Lip}. Assume that $N>K>4|\cO|$ and let $\Delta=1/4-|\cO|/K$. Then, with probability larger than $1-2\exp\left(-2\Delta^2K\right),$ we have
$$R(\hat f_{\text{MOM},K})\leq \inf_{f\in F}R(f) +  4L\max\pa{\frac{4\cR(F)}{ N}, 2\param{2}\sqrt{\frac{K}{ N}}}\enspace.$$
\end{Th}

Theorem~\ref{thm:MOM} is proved in Section~\ref{sec:ProofThm1}. Compared to Theorem~\ref{Lem:ERMOutliers}, $\hat f_{\text{MOM},K}$ achieves the same  rate $(\cR(F)/ N)\vee(\sqrt{K/ N})$ under the same conditions on the number of outliers with the same exponential control of the probability as for the ERM $f_{\text{ERM}}^{0-1}$. The main difference is that the loss function may be unbounded, which is often the case in practice. Moreover, unlike classical analysis of ERM obtained by minimizing an empirical risk associated with a convex surrogate loss function, we only need a second moment assumption on the class $F$.

These theoretical improvements have already been noticed in previous works \cite{MR3378468,subgaussian,MOM1,LugosiMendelson2017,LugosiMendelson2017-2,LugosiMendelson2016,shahar_mom_1,MOM2}. 
Contrary to tournaments of \cite{LugosiMendelson2016}, Le Cam MOM estimators of \cite{MOM1} or minmax MOM estimators \cite{MOM2}, Theorem~\ref{thm:MOM} does not require the small ball assumption on $F$ but only shows ``slow rates" of convergence.
These slow rates are minimax optimal in the absence of a margin or Bernstein assumption \cite{MR2240689,MR1765618}.
Proof of Theorem~\ref{thm:MOM} does not enable fast rates to be obtained.  
Indeed, the non-linearity of the median excludes the possibility of using localization techniques leading to these fast rates.
However, we show in the simulation study (cf. left side picture of Figure~\ref{fig:rates}) that fast rates seem to be reached by the MOM minimizer.

The main advantage of our approach is its simplicity, we just have to replace empirical means by their MOM alternative in the definition of the ERM. Moreover, as expected, this simple alternative to ERM yields a systematically way to modify algorithms designed for approximating the ERM. The resulting ``MOM version'' of these algorithms are more robustness and faster than their original ``ERM version''. Before illustrating these facts on simulations, let us describe algorithms approximating MOM minimizers. 

\section{Computation of MOM minimizers}\label{sec:Comp}
In this section, we present first a generic algorithm to provide a MOM version of descent algorithms minimizing the empirical risk. We show how classical convergence proofs can be adapted to these algorithms.

\subsection{MOM algorithms}
The general idea is that any descent algorithms such as gradient descent, Newton method, alternate gradient descent, etc. (cf. \cite{moulines2011non,BubeckOptim,MR2061575,MR3025128}) can easily be turned into a robust MOM-version. To illustrate this idea, a basic gradient descent is analyzed in the sequel. 

The choice of blocks greatly influences the practical performance of the algorithm. 
In particular, a recurring flaw is that iterations tend to get stuck in local minima, which greatly slows the convergence of the alogorithm.
To overcome this default and improve the stability of the procedure, we modify the partition at each iteration by drawing it uniformly at random, cf. step \textbf{2} of Algorithm~\ref{fig:mom_algo}.

Let $\cS_N$ denote the set of permutations of $\{1,\ldots,N\}$. For each $\sigma\in \cS_N$, let $B_0(\sigma)\cup\cdots\cup B_{K-1}(\sigma)=\{1, \ldots, N\}$ denote an equipartition of $\{1, \ldots, N\}$ defined for all $ j \in \llb 0,K-1\rrb$ by
$$B_j(\sigma) =\{\sigma(Kj+1),\dots,\sigma(K(j+1))\} = \sigma \left(\{Kj+1, \ldots, K(j+1)\}\right)\enspace.$$

To simplify the presentation, let us assume the class $F$ to be parametrized $F=\{ f_u : u \in \R^p\}$, for some $p\in \mathbb{N}^*$. Let's assume that the function $u\mapsto f_u$ is as regular as needed and convex (a typical example is $f_u(x) = \inr{u, x}$ for all $x\in\R^p$). Denote by $\nabla_u\ell_{f_u}$ the gradient or subgradient of $u\mapsto f_u$ in $u\in\R^p$. The step-sizes sequence is denoted by  $(\eta_t)_{t\ge 0}$ and satisfies the classical conditions:  $\sum_{t=1}^{\infty}\eta_t = \infty$ and $\sum_{t=1}^{\infty}\eta_t^2 < \infty$. Iterations will go on until a stopping time $T\in \mathbb{N}^*$ has been achieved. With these notations, a generic MOM version of a gradient descent algorithm (with random choice of blocks) is detailed in Algorithm~\ref{fig:mom_algo}.

\vspace{0.6cm}
\begin{center}
\begin{algorithm}[H]\label{fig:mom_algo}
\SetKwInOut{Input}{input}\SetKwInOut{Output}{output}\SetKw{Or}{or}\SetKw{In}{in}
\SetKw{Return}{Return}
\Input{$u_0\in \R^p$, $K\in \llb 3, N/2 \rrb$, $T \in \mathbb{N}^*$ and $(\eta_t)_{t\in \{0,\ldots,T-1\}}\in\R^{T}_+$}
\Output{approximated solution to the minimization problem \eqref{eq:defMOMK}}  
\BlankLine

\For{$t =0 ,\cdots, T-1$}{
	choose a permutation at random: $\sigma_t \sim \text{Unif}(\mathcal{S}_N)$,\\
	build a partition of the dataset: $B_0(\sigma_t),\dots,B_{K-1}(\sigma_t)$,\\
	find a median block: $k_{med}(t)$ s.t. $\MOM{K}{\ell_{f_{u_t}}}=P_{B_{k_{med}(t)}(\sigma_t)}\left(\ell_{f_{u_t}}\right)$,\\
	do a descent step
$$u_{t+1}=u_{t}-\eta_t \sum_{i \in B_{k_{med}(t)}(\sigma_t)}\nabla_{u_t}\ell_{f_{u_t}}(X_i,Y_i).$$
}
 \Return $u_T$
 \caption{MOM gradient descent algorithm.}
\end{algorithm} 
\end{center}

\begin{Remark}[MOM gradient descent algorithm and stochastic block gradient descent]
Algorithm~\ref{fig:mom_algo} can be seen as a stochastic block gradient descent (SBGD) algorithm minimizing the function $t\to \E \ell_t(X, Y)$ using a given dataset. The main difference with the classical SBGD is that the choice of the block along which the gradient direction is performed is chosen according to a centrality measure computed thanks to the median operator in step \textbf{4} of Algorithm~\ref{fig:mom_algo}. 
\end{Remark}
In Section~\ref{simu}, we use MOM ideas (as in the generic Algorithm~\ref{fig:mom_algo}) to construct MOM versions for various classical algorithms such as Perceptron, Logistic Regression, Kernel Logistic Regression, SGD Classifiers or Multi-layer Perceptron. 



\subsection{Convergence properties}\label{sec:conv_MOM_algo}
We denote by $\cD_N=\{(X_1,Y_1),\dots,(X_N,Y_N)\}$ the dataset. We consider the following set of assumptions: 
\begin{Assumption}There exist $L>0$ and a sequence of decreasing positive numbers $(\eta_t)_{t\geqslant 1}$, such that
 \begin{enumerate}[label=\textbf{H-\arabic*}]
\item 
for all $u \in \R^p$,and $P$-almost all $x,y\in \cX \times\{-1, 1\}$, $\norm{\nabla_u\ell_{f_u}(x,y)}_2\le L$, \label{hyp:lipshitz}
\item for  almost all datasets $\cD_N=\{(X_i,Y_i):i=1, \ldots, N\}$, there exists a unique minimum $u_{min} =\argmin_{u \in \R^p }\E_\sigma[\MOM{K}{\ell_{f_u}}| \cD_N]$ where $\E_\sigma$ denotes the expectation with respect to $\sigma$, the uniformly distributed permutation used to construct the blocks of data,\label{hyp:existance}
\item $\sum_{t=1}^{\infty}\eta_t^2 <\infty$,\label{hyp:seriesconv}
\item $\sum_{t=1}^{\infty}\eta_t=\infty$,\label{hyp:seriesdiv}
\item for  almost all datasets $\cD_N=\{(X_i,Y_i):i=1, \ldots, N\}$, for Lebesgue almost all $u \in \R^p$, for all $\varepsilon>0$, $\inf_{\norm{u_{min}-u}_2>\varepsilon}(u_{min}-u)^T \E_\sigma[\nabla_u\MOM{K}{\ell_{f_u}}|\cD_N]> 0.$ \label{hyp:convex}
\item for  almost all datasets $\cD_N=\{(X_i,Y_i):i=1, \ldots, N\}$ and Lebesgue almost all $u \in \R^p$, there exists an open convex set $B$ containing $u$ such that for any equipartition of $\{1, \ldots, N\}$ into $K$ blocks $B_1, \ldots, B_K$ there exists  $k_{med}\in \{1,\cdots, K \}$ such that for all $v\in B$, $P_{B_{k_{med}}}(\ell_{f_{v}})\in \MOM{K}{\ell_{f_{v}}}$. \label{hyp:local_conv}
\end{enumerate}
\end{Assumption}

Recall that the dataset $\cD_N=\{(X_i,Y_i):i=1, \ldots, N\}$ may contain outlier, which might be deterministic. We use the terminology ``for almost all datasets $\cD_N=\{(X_i,Y_i):i=1, \ldots, N\}$ such assumption holds" when the latter holds for $P^{|\cI|}$-almost all informative data $(X_i, Y_i)_{i\in \cI}$ and for every outliers $(X_i, Y_i)_{i\in\cO}\in (\cX\times \cY)^{|\cO|}$.

\begin{Th}[Convergence of the algorithm]
\label{thm:convergence}
Grant Assumptions \ref{hyp:lipshitz} to \ref{hyp:local_conv}, then the MOM gradient descent algorithm~\ref{fig:mom_algo} $(u_t)_{t\in\bN}$ converges: for almost all $\cD_N$, we have
$$\norm{u_T-u_{min}} \xrightarrow[T\to \infty]{a.s}0 $$
(where the a.s. convergence is with respect to the randomness on the choice of data partitions at each step).
\end{Th}
Note that Theorem~\ref{thm:convergence} insures the convergence of the MOM version of the gradient descent algorithm with random blocks as defined in Algorithm~\ref{fig:mom_algo} toward a minimizer of the conditional expectation $u\to\E_\sigma[\MOM{K}{\ell_{f_u}}| \cD_N]$. It does not prove its convergence toward a minimizer of the MOM criteria $u\to \MOM{K}{\ell_{f_u}}$. Nevertheless,  we can expect a strong concentration of $\MOM{K}{\ell_{f_u}}$ around $\E_\sigma[\MOM{K}{\ell_{f_u}}| \cD_N]$ by concentration of random permutations, see \cite{mcdiarmid_2002}.

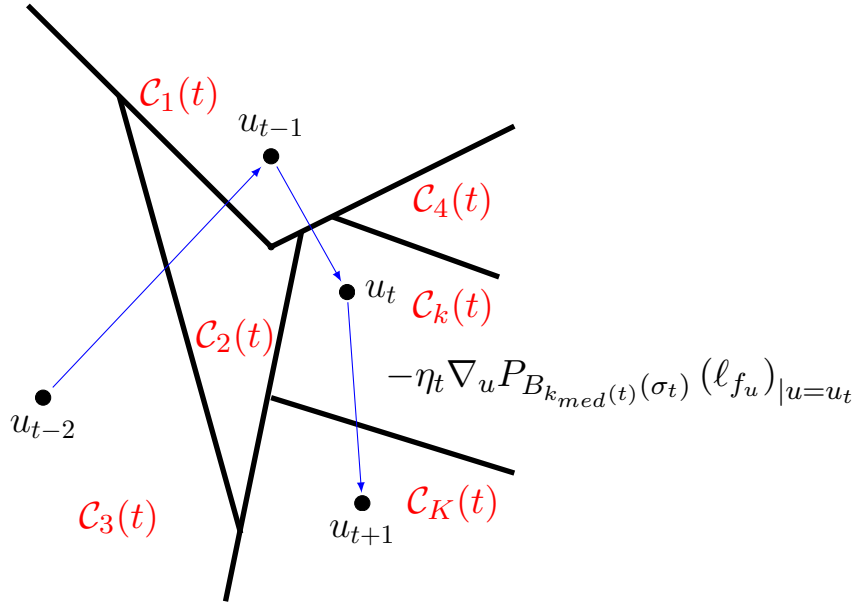
\begin{figure}[h]
\centering
\centering
\begin{tikzpicture}[scale=0.2]
\filldraw (0,0) circle (0.2cm);
\draw[line width = 2pt] (-0.1,0.1) -- (-16, 16);
\draw[line width = 2pt] (0,0) -- (16,8);
\draw[line width = 2pt] (-10,10) -- (-2,-19);
\draw[line width = 2pt] (2,1) -- (-3,-23.5);
\draw[line width = 2pt] (0,-10) -- (16,-15);
\draw[line width = 2pt] (4,2) -- (15,-2);

\filldraw (-15,-10) circle (5mm) node (utm2) {} node[anchor=north,fill=white,yshift=-0.1cm]{\Large $u_{t-2}$};
\filldraw (0,6) circle (5mm) node (utm1) {} node[anchor=south,fill=white,yshift=0.1cm]{\Large $u_{t-1}$};
\filldraw (5,-3) circle (5mm) node (ut) {} node[anchor=west,fill=white,xshift=0.1cm]{\Large $u_{t}$};
\filldraw (6,-17) circle (5mm) node (utp1) {} node[anchor=north,fill=white,yshift=-0.1cm]{\Large $u_{t+1}$};

\draw[-latex, blue] (utm2) -- (utm1);
\draw[-latex,blue](utm1) -- (ut);
\draw[-latex,blue] (ut) -- (utp1);
\draw (23,-9) node {\scalebox{1.3}{$-\eta_t \nabla_{u}P_{B_{k_{med}(t)}(\sigma_t)}\left(\ell_{f_u}\right)_{|u=u_t}$}};
\draw[red] (-6,10) node {\scalebox{1.3}{$\cC_1(t)$}};
\draw[red] (-2.3,-6) node {\scalebox{1.3}{$\cC_2(t)$}};
\draw[red] (-10,-18) node {\scalebox{1.3}{$\cC_3(t)$}};
\draw[red] (12,3) node {\scalebox{1.3}{$\cC_4(t)$}};
\draw[red] (12,-4) node {\scalebox{1.3}{$\cC_k(t)$}};
\draw[red] (12,-17) node {\scalebox{1.3}{$\cC_K(t)$}};
\end{tikzpicture}
\caption{Partition of $\R^p$ at step $t$ by the median operator and iteration number $t-2, t-1, t$ and $t+1$ of the MOM gradient descent algorithm. Under Assumption~\ref{hyp:local_conv}, there is a natural descent direction given at step $t$ by $-\nabla_{u}(u\to P_{B_{k_{med}(t)}(\sigma_t)}\left(\ell_{f_u}\right)))_{|u=u_t}$. 
}
\label{fig:MOM_descent}
\end{figure}

\begin{Remark}\label{rem:sgd}[MOM gradient descent, SGD and assumption~\ref{hyp:local_conv}]
All assumptions in Theorem~\ref{thm:convergence} are the classical assumption used to prove the convergence of the classical SGD except for \ref{hyp:local_conv} which is specific to MOM algorithms. 

Indeed, the picture of the MOM gradient descent algorithm is pretty simple and depicted in Figure~\ref{fig:MOM_descent}. At every step $t$, the median operator makes a partition of $\R^p$ into $K$ cells $\cC_k(t) = \{u\in\R^p: \MOM{K}{\ell_{f_u}} = P_{B_k}\ell_{f_u}\}$ for $k=1, \ldots, K$ -- this partition changes at every step because the blocks $B_1, \ldots, B_K$ are chosen randomly at the beginning of every step according to the random partition $\sigma_t$. We want every iteration $u_t$ of the MOM algorithm to be in the interior of a cell and not on a frontier in order to differentiate the objective function $u\to \MOM{K}{\ell_{f_u}}$ at $u_t$. This is indeed the case under \ref{hyp:local_conv}, given that in that case, there is an open neighbor $B$ of $u_t$ such that for all $v\in B$, $\MOM{K}{\ell_{f_v}} = P_{B_k}\ell_{f_v}$ where the index $k=k_{med}$ of the block is common to every $v\in B$. Therefore, to differentiate the objective function $u\to \MOM{K}{\ell_{f_u}}$ at $u_t$ one just needs to differentiate $u\to P_{B_k}\ell_{f_u}$ at $u_t$.  
The objective function to minimize is differentiable almost everywhere under \ref{hyp:local_conv} and a gradient of the objective function is given by $\nabla(u\to P_{B_k}\ell_{f_u})_{|u=u_t}$.
\end{Remark}

\begin{figure}[h]
\begin{center}
\includegraphics[scale=0.7]{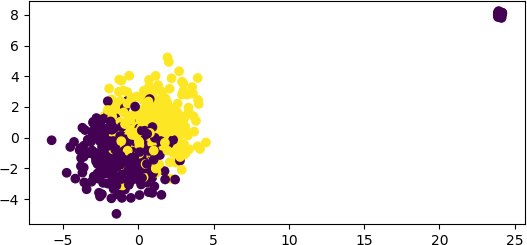}
\caption{Scatter plot of 630 samples from the training dataset (600 informative data, 30 outliers), the color of the points correspond to their labels.\label{fig:dataset}}
\end{center}
\end{figure}

\begin{Remark}[Assumption~\ref{hyp:convex}]
The assumption~\ref{hyp:convex} illustrate why we used a random shuffle of the data at each step of the algorithm, the principle is that the MOM risk in general may present local minima and hence the algorithm without the shuffle would stop in a local minimum. By using the random shuffle at each step, it seems that the problem has a unique minimum, at least experiment seems to say so. For example, in figure~\ref{fig:loss} we see the plot of the loss function when we change one of the parameters (the intercept) of a linear classifier for the classification problem depicted in figure~\ref{fig:dataset}.

\begin{figure}[h]
\begin{center}

    \subfloat[Plot of MOM risk]{
      \includegraphics[width=0.45\textwidth]{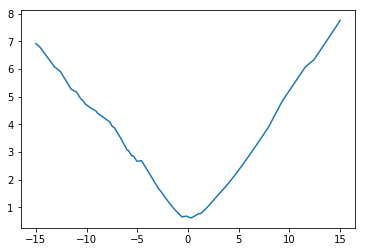}
                         }
    \subfloat[Plot of expected MOM risk]{
      \includegraphics[width=0.45\textwidth]{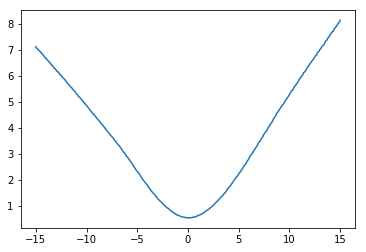}
                         }

\caption{Plot of the two loss functions the MOM risk $\MOM{K}{\ell_{f_b}}$ and the expected MOM risk $\E[\MOM{K}{\ell_{f_b}}|\cD_N]$ (estimated using Monte Carlo method with $300$ samples), as a function of $b$ where $\ell$ is the logistic loss and $f_b$ is the linear classifier $f_b(x)=bx_1+x_2+x_3, \forall x=(x_j)_{j=1}^3\in\R^3$ and $b\in \R$. \label{fig:loss}}
\end{center}
\end{figure}

Remark that assumption~\ref{hyp:convex} is not convexity, and indeed in figure~\ref{fig:loss} the expected loss is not convex but it verifies assumption~\ref{hyp:convex} nonetheless.
\end{Remark}

\begin{Example}
Let us give a relatively general example where Assumption~\ref{hyp:local_conv} is satisfied. Let $B_1\cup\cdots\cup B_K=\{1, \ldots, N\}$ be an equipartition and let $\psi$ be defined for all $x=(x_i)_{i=1}^N\in \mathbb{R}^N$ by,
$$\psi_p(x)=\MOM{K}{f_p(x)}=\text{median}\left( \frac{K}{N}\sum_{i\in B_k}f_p(x_i), 1\le k\le K \right)= P_{B(K/2)(p)}(f_p),$$where for all blocks $B\subset \{1, \ldots, N\}, P_Bf_p = |B|^{-1}\sum_{i\in B} f_p(x_i)$ and the blocks $B(k)(p),k=1, \ldots, K$ are rearranged blocks defined such that $P_{B(1)(p)}(f_p)\geq\dots\geq P_{B(K)(p)}(f_p)$. Proposition~\ref{prop:example_gradient} shows that Assumption~\ref{hyp:local_conv} is satisfied in several situations.

Before stating the proposition, recall that the hazard function of an absolutely continuous real-valued random variable $Z$ is defined for all $x\in\R$ by $h(x)=f(x)/(1-F(x))$ where $f$ is the density function of $Z$ and $F$ is its cumulative distribution function. 

\begin{Prop}
\label{prop:example_gradient}Let $X_1, \ldots, X_N$ be $N$ real-valued random variables and let $(f_p)_{p \in \R^d}$ be a family of functions with values in $\R$ , if for all $x\in\R$, the function $p\mapsto f_p(x)$ is Lipshitz and the probability distribution of a mean of a block $P_{B(k)(p)}(f_p)$ 
 is a random variable (recall that the choice of block is random) with a bounded  hazard function, then with probability $1$, Assumption~\ref{hyp:local_conv} is satisfied, in particular, the partial derivative of $p\mapsto \psi_{f_p}((X_i)_{i=1}^N)=MOM_K(f_p((X_i)_{i=1}^N))$ with respect to the $j^{th}$ coordinate is given for almost all $X_1, \ldots, X_N$ by
$$\partial_j\psi_{f_p}((X_i)_{i=1}^N)= \frac{K}{N}\sum_{i\in B(K/2)(p)}\partial_jf_p(X_i)$$
where $\partial_j$ denote the derivative with respect to the $j^{th}$ coordinate of $p$. 
\end{Prop}

Proposition~\ref{prop:example_gradient} is proved in Section \ref{sec:Proof}.
The function $h$ is bounded for a large panel of probability distributions. For example, an exponential variable with parameter $\lambda$ has a constant hazard function equal to $\lambda$; a Cauchy random variable has also a bounded hazard function. Intuitively, this property of derivation of the median of means is more easily verified when the distribution of $P_{B(k)(p)}(f)$ is heavy tailed (i.e. the mean blocks are far from one another) and in fact we can verify that the hazard function of the Gaussian is unbounded and hence this proposition cannot be used on Gaussian random variables.

\end{Example}

\subsection{Complexity of MOM risk minimization algorithms}
\label{complexity}
Let $C(m)$ be the computational complexity of a single ERM gradient descent update step on a dataset of size $m$ and let $L(m)$ be the computational complexity of the evaluation of the loss function on a dataset  containing $m$ data. Here the computational complexity is simply the number of basic operations needed to perform a task \cite{MR2500087}.

For each epoch, we begin by computing the ``MOM empirical risk''. We perform $K$ times $N/K$ evaluations of the loss function, then we sort the $K$ means of these blocks of loss to finally get the median. The complexity of this step is then $O(KL(N/K)+K\ln(K))$, assuming that the sort algorithm is in $O(K\ln(K))$ (like \textit{quick sort} \cite{hoare1962quicksort}). Then we do the ERM gradient step on a sample of size $N/K$. Hence, the time complexity of this algorithm is 
$$O(T(KL(N/K)+K\ln(K)+C(N/K))).$$
\begin{Example}[Linear complexity ERM algorithms]
For example, if the ERM gradient step and the loss function evaluation have linear complexity -- like Perceptron or Logistic Regression -- the complexity of the MOM algorithm is $O(T(N+K\log(K)))$ against $O(TN)$ for the ERM algorithm. Therefore, the two complexities are of the same order and the only advantage of MOM algorithms lies in their robustness to outliers.
\end{Example} 
\begin{Example}[Super-linear complexity ERM algorithms]
If on the other hand the complexity is more than linear as for Kernel Logistic Regression (KLR), taking into account the matrix multiplications whose complexity can be found in \cite{DBLP:journals/corr/Gall14}, the complexity of the MOM version of KLR, due to the additional need of the computation of the kernel matrix, is
 $O(N^2+T(N^2/K+K\log(K)+(N/K)^{2.373}))$ against  $O(TN^{2.373})$ for the ERM version. MOM versions of KLR are therefore faster than the classical version of KLR on top of being more robust. This advantage comes from the fact that MOM algorithms work on blocks of data instead on the entire dataset at every step. More informations about Kernel Logistic Regression can be found in \cite{Roth01probabilisticdiscriminative} for example.
\end{Example}

In this last example, the complexity comes in part from the evaluation of the kernel matrix that can be computationally expensive. Following the idea that MOM algorithms are performing ERM algorithm restricted to a wisely chosen block of data, then one can modify our generic strategy in this particular example to reduce drastically its complexity. The idea here is that we only need to construct the kernel matrix on the median block. The resulting algorithm, called Fast KLR MOM is described in Figure~\ref{fig:mom_algo_2}. 

In Figure~\ref{fig:mom_algo_2}, we compute only the block kernel matrices, denoted by $N^1,\dots,N^k$ and constructed from the samples in the block $B_k$. We also denote by $N_i^k$ the $i^{th}$ row in $N^k$.

\begin{center}
\begin{algorithm}[H]\label{fig:mom_algo_2}
\SetKwInOut{Input}{input}\SetKwInOut{Output}{output}\SetKw{Or}{or}\SetKw{In}{in}
\SetKw{Return}{Return}
\Input{$\alpha_0\in \R^p$, $K\in \llb 3, N/2 \rrb$, $T \in \mathbb{N}^*$, $(\eta_t)_{t\in \{0,\ldots,T-1\}}\in\R^{T}_+$, $\beta \in \R_+^*$, $\kappa:\cX\times \cX \to \R$ a positive definite kernel and a bloc decomposition $B_1,\dots,B_K$ of $\{1,\dots,N\}$.}
\Output{approximated solution to the minimization problem \eqref{eq:defMOMK} among the KLR classifiers}  
\BlankLine
Construct the bloc Kernel matrices $N^k=(\kappa(X_i,X_j))_{i,j \in B_k}$ for $1\le k\le K$,\\
\For{$t =0 ,\cdots, T-1$}{
	find a median block: $k_{med}(t)$ s.t.  $\MOM{K}{\ell_{f_{\alpha_t}}}=P_{B_{k_{med}(t)}}\left(\ell_{f_{\alpha_t}}\right)$ with\\
	$$P_{B_k}(\ell_{f_{\alpha_t}}) =\frac{1}{|B_k|}\sum_{i\in B_k}\ln(1+e^{-N_i^k\alpha_t^k Y_i})+\beta \sum_{k=1}^K (\alpha_t^k)^T N^k \alpha_t^k,$$
	where $\alpha_t^k$ is the vector in $\R^{|B_k|}$ made of the coordinates of $\alpha_t$ with indices in $B_k$.\\ 
	Do an IRLS descent step for KLR with weight matrice $W_{k_{med}(t)}$, design matrice $X_{k_{med}(t)}$ and labels $y_{k_{med}(t)}$ on $B_{k_{med}(t)}$
	\begin{align*}\alpha_{t+1}^{k_{med}(t)}&=\alpha_{t}^{k_{med}(t)}(1-\eta_t)+\eta_t (X_{k_{med}(t)}^TW_{k_{med}(t)}X_{k_{med}(t)})^{-1}X_{k_{med}(t)}^TW_{k_{med}(t)}y_{k_{med}(t)}.\\
	\alpha_{t+1}^{k}&=\alpha_{t}^{k}(1-\eta_t), \quad \forall k \neq k_{med}(t).
	\end{align*}
}
 \Return $\alpha_T,N_{k_{med}(T)}$
 \caption{Description of Fast KLR MOM algorithm.}
\end{algorithm} 
\end{center}

%

There are several drawbacks in the approach of Algorithm~\ref{fig:mom_algo_2}. First, the blocks are fixed at the beginning of the algorithm; therefore the algorithm needs a bigger dataset to work well and it may converge to a local minimum. Nonetheless, from the complexity point of view, this algorithm will be much faster than both the classical KLR and MOM KLR (see below for a computation of its complexity) 
which is important given the growing use of kernel methods on very large databases for example in biology. The choice of $K$ should ultimately realize a trade-off between complexity and performance (in term of accuracy for example) 
when dealing with big databases containing few outliers.
\begin{Example}[Complexity of Fast KLR-MOM algorithm]
The complexity of Fast KLR-MOM is $O(N^2/K+T(N^2/K+K\log(K)+(N/K)^{2.373}))$ against $O(TN^{2.373})$ for the ERM version.
\end{Example}
\section{Implementation and Simulations}
\label{simu}
\subsection{Basic results on a toy dataset}
The toy model we consider models outliers due to human or machine errors we would like to ignore in our learning process. It is also a dataset corrupted to make linear classifiers fail. The dataset is a 2D dataset constituted of three ``labeled Gaussian distribution''. Two informative Gaussians $\mathcal{N}((-1,-1),1.4I_2)$ and $\mathcal{N}((1,1),1.4I_2)$ with label respectively $1$ and $-1$ and one outliers Gaussian $\mathcal{N}((24,8),0.1I_2)$ with label $1$. In other words, the distribution of informative data is given by $\cL(X|Y=1) = \cN((-1,-1),1.4I_2)$, $\cL(X|Y=-1) = \cN((1,1),1.4I_2)$ and $\bP(Y=1)=\bP(Y=-1)=1/2$. Outliers data have distribution given by $Y=1$ a.s. and $X\sim \mathcal{N}((24,8),0.1I_2)$.

The algorithms we study are the MOM adaptations of Perceptron, Logistic Regression and Kernel Logistic Regression. 

Based on our theoretical results, we know that the number of blocks $K$ has to be larger than $4$ times the number of outliers for our procedure to be on the safe side. The value $K=120$ is therefore used in all subsequent applications of MOM algorithms on the toy dataset except when told otherwise. To quantify performance, we compute the miss-classification error on a clean dataset made of data distributed like the informative data. 

\begin{figure}[h!]
\begin{center}
\includegraphics[scale=0.7]{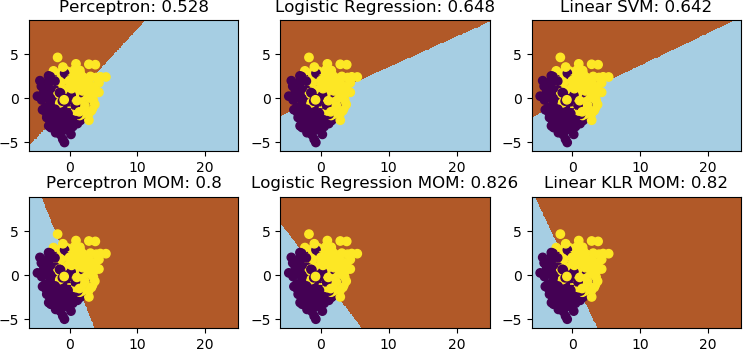}
\caption{Scatter plot of 500 samples from the test dataset (500 informative data), the color of the points correspond to their labels and the background color correspond to the prediction. The score in the title of each subfigure is the accuracy of the algorithm.\label{fig:classif_comparaison}}
\end{center}
\end{figure}

For Kernel Logistic Regression, we study here a linear kernel because outliers in this dataset are clearly adversarial when dealing with linear classifiers. The algorithm can also use more sophisticated kernels, a comparison of the MOM algorithms with similar ERM algorithms is represented in figure \ref{fig:classif_comparaison}, the ERM algorithms are taken from the python library scikit-learn \cite{sklearn} with their default parameters.

Figure \ref{fig:classif_comparaison} illustrates resistance to outliers of MOM's algorithms compared to their classical version.

\begin{figure}[h]
  \begin{center}

\includegraphics[scale=0.7]{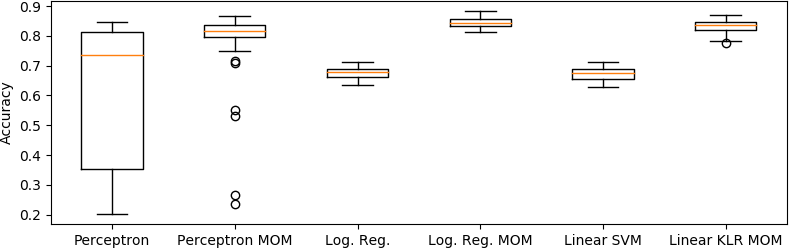}
 
\caption{Comparison of the MOM algorithms  and their counterpart with the boxplots of the accuracy on the test dataset from 50 runs of the algorithms on 50 sample of the training/test toy dataset (one run for each dataset sampled)\label{fig:boxplot1}.}
  \end{center}
\end{figure}

\noindent
These first results are completed in Figure \ref{fig:boxplot1} where we computed accuracy on several run of the algorithms. These results confirm the visual impression of our first experiment.

Finally, we illustrate our results regarding complexities of the algorithms on a simulated example. MOM algorithms have been computed together with state-of-the art algorithms from scikit-learn \cite{sklearn} (we use Random forest, SVM classifier as well as SGD classifier optimizing Huber loss which entail a robustness in $Y$ but not in $X$, see \cite[Chapter~7]{robuststat}) on a simulated dataset composed of two Gaussian blobs  $\mathcal{N}((-1,-1),1.4I_2)$ and $\mathcal{N}((1,1),1.4I_2)$ with label respectively $1$ and $-1$. We sample $20 000$ points for the training dataset and $20 000$ for the test dataset. The parameters used in the algorithms are those for which we obtained the optimal accuracy, (this accuracy is illustrated in the next section).
 Time of training plus time of evaluation on the test dataset are gathered in Figure~\ref{fig:complexities_table}. 

\begin{figure}[h]
\begin{center}
\begin{tabular}{|c|c|c|c|c|}
  \hline
   Algorithm&Perceptron MOM &  Log. Reg. MOM & KLR MOM & Fast KLR MOM \\
  \hline
  Time (s)&$1.06 $ & $ 1.05$& $ 13.6$&$1.2 $  \\
  \hline
  \hline 
 Algorithm & Rand. Forest &  SVM &  SGD Hub. loss&\\
 \hline
   Time (s)&$ 0.21$ & $ 9.0$&$0.0078 $ & \\
  \hline
\end{tabular}
\caption{Time of different algorithms on a simulated dataset \label{fig:complexities_table}.}
\end{center}
\end{figure}

Not surprisingly, very efficient versions of linear algorithms from Python's library are extremely fast (results are sometimes provided before we even charged the dataset in some experiments). The performance of our algorithm are nevertheless acceptable in general (around 5 times longer than random forest for example). The important fact here is that non linear algorithms such as SVM take much more time to provide a result. FAST KLR MOM is able to reduce substantially the execution time of SVM with comparable predictive performance.  

\subsection{Applications on real datasets}
We used the HTRU2 dataset, also studied in \cite{pulsar}, that is provided by the UCI Machine Learning Repository. The goal is to detect pulsars (a rare type of Neutron star) based on radio emission detectable on earth from which features are extracted to gives us this dataset. The problem is that most of the signal comes from noise and not pulsar, the goal is then to classify pulsar against noise, using the $17$ $898$ points in the dataset.

\begin{figure}[h]
  \begin{center}

      \includegraphics[scale=0.5]{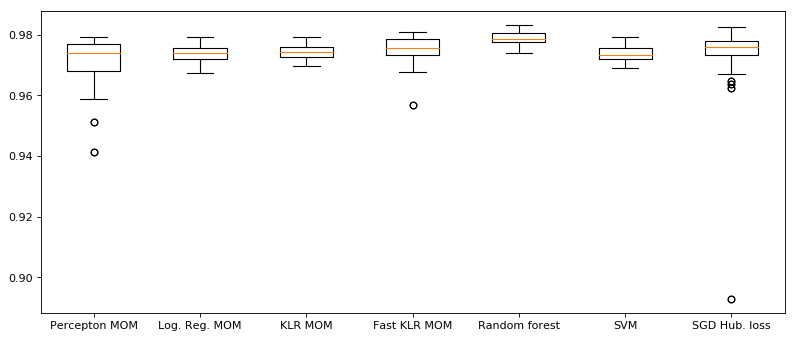}
 
    \caption{Comparison of the MOM algorithms and common algorithms with the boxplots and the medians of the accuracy $\frac{1}{n}\sum_{i=1}^n \1\{\hat f(x_i)=y_i\} $ on the test dataset from 50 runs of the algorithms on 100 sample of a $4/5$ cut of the dataset HTRU2 (one run is trained on a sample of $4/5$ of the dataset and tested on the remaining $1/5$)\label{fig:boxplot_pulsar}}
  \end{center}
\end{figure}

The accuracy of different algorithms is obtained using on several runs of the algorithms each using 4/5 of the datasets for training and $1/5$ for testing algorithms. Boxplots presenting performance of various algorithms are displayed in Figure \ref{fig:boxplot_pulsar}. To improve performance, RBF kernel was used both for KLR MOM and Fast KLR MOM.

\subsection{Outlier detection with MOM algorithms}
When we run MOM version of a descent algorithm, we select at each step a block of data points realizing the median of a set of ``local/block empirical risk'' at the current iteration of the algorithm. The number of times a point is selected by the algorithm can be used as a depth function measuring reliability of the data. Note that this definition of depth of a data point has the advantage of taking into account the learning task we want to solve, that is the loss $\ell$ and the class $F$. It means that outliers are considered w.r.t. the problem we want to solve and not w.r.t. some a priori notion of centrality of points in $\R^d$ unrelated with the problem considered at the beginning.

\begin{figure}[h]
  \begin{center}

      \includegraphics[scale=0.7]{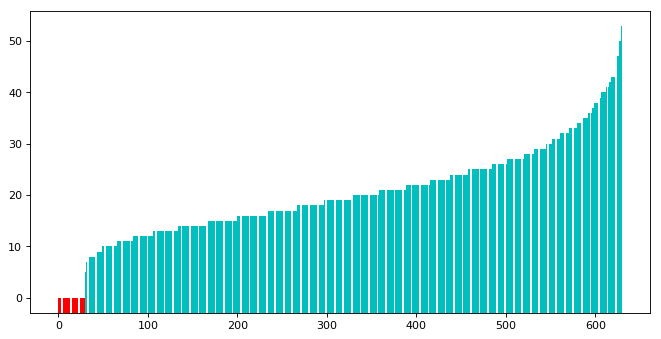}
 
    \caption{Sorted Histogram of the score (number of times a data belongs to the selected median block) of each points in a Logistic Regression MOM algorithm on a toy dataset. Red is an outlier and blue is an informative sample. $K=120$ and $T=2000$. \label{fig:outlier1}}
  \end{center}
\end{figure}

We apply this idea on the toy dataset with the Logistic Regression MOM algorithm. Results are gathered in a sorted histogram given in Figure \ref{fig:outlier1}. Red bars represent outliers in the original datasets. 

Quite remarkably, outliers are in fact those data that have been used the smallest number of times. 
The method targets a very specific type of outliers, those disturbing the classification task at hand. If there was a point very far away from the bulk of data but in the half-space of its label, it wouldn't be detected. 

This detection algorithm doesn't scale well when the dataset gets bigger as a large number of iterations is necessary to choose each point a fair number of times. For bigger datasets, we suggest to adapt usual outlier detection algorithms \cite{outlier_analysis}. We emphasize that clustering techniques and K-Means are rather easy to adapt in a MOM algorithm and detect points far from the bulk of data. This technique might greatly improve usual K-Means as MOM K-Means is more robust.

\begin{figure}[h]
  \begin{center}

     \includegraphics[scale=0.42]{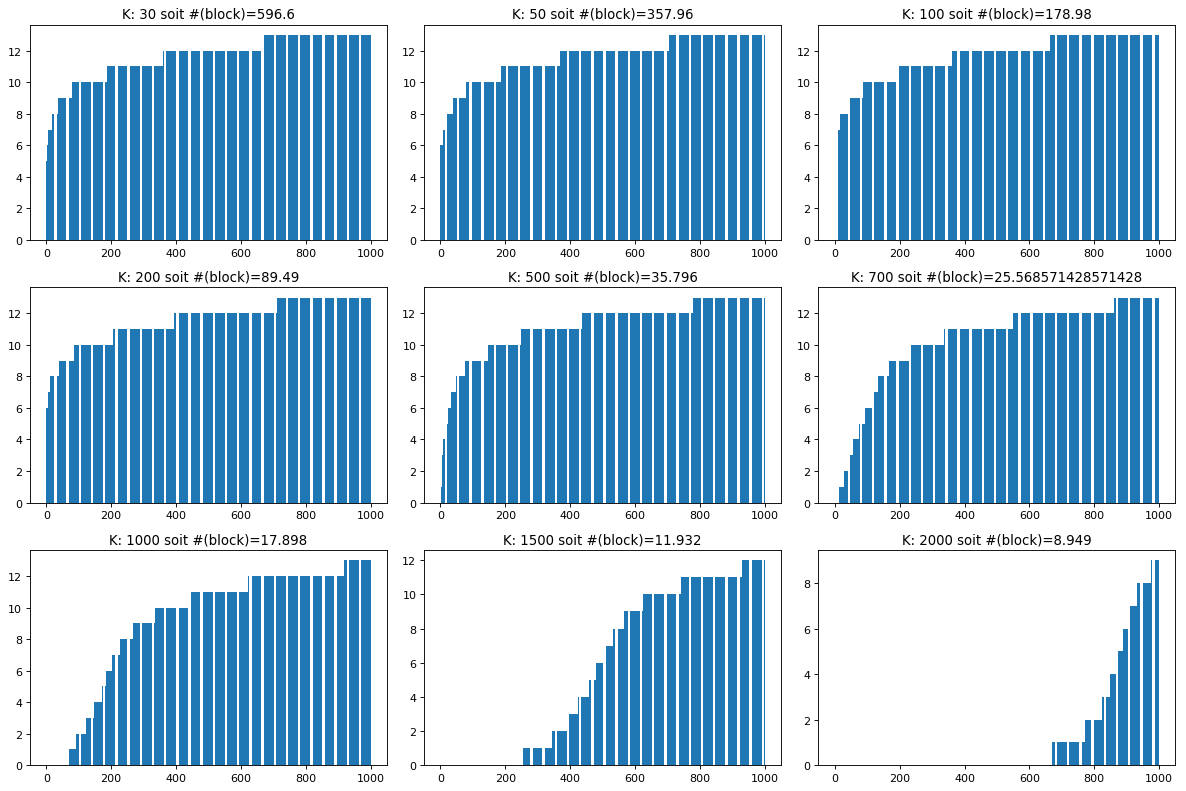}
 
    \caption{Sorted Histogram on the score (number of times a data is selected in a median block) of each points in a Logistic Regression MOM algorithm on the pulsar dataset for various values of $K$ and $T=20 \times K$  (only the 1000 smaller counts among the 17898 sample of the pulsar dataset are represented). \label{fig:outlier_pulsar}}
  \end{center}
\end{figure}

Let us now analyze the effect of $K$ on the outlier detection task. The histogram of the $1000$ smaller counts of points of HTRU2 dataset as $K$ gets bigger is plotted in Figure \ref{fig:outlier_pulsar}.

It appears from Figure \ref{fig:outlier_pulsar} that $K$ measures the sensitivity of the algorithm. Severe outliers (as in the toy example) are detected for small $K$ while mild outliers are only discovered as $K$ gets bigger.

It seems therefore that the optimal choice of $K$ in MOM depends on the task one is interested in. For classification, $K$ should be as small as possible to get better risk bounds (but it still should be larger than the number of outliers) whereas for detecting outliers we may want to choose $K$ much larger to even detect an outlier, (but it should also be small enough for the underlying classification to perform correctly). As a proof of concept, for Pulsar database, we got optimal results choosing $K=10$  for classification whereas we only detect a significant amount of outliers when $K$ is around $1000$.

\section{Proofs}\label{sec:Proof}

\subsection{Proof of Theorem~\ref{Lem:ERMOutliers}}\label{sec:ProofLem1}
We adapt Vapnik's classical analysis \cite{MR1641250} of excess risk bound of ERM to a dataset corrupted by outliers. We first recall that $f^*\in\argmin_{f\in\cF} R^{0-1}(f)$ and for all $f\in\cF$, the excess loss function of $f$ is $\cL_f^{0-1}=\ell^{0-1}_f-\ell^{0-1}_{f^*}$. For simplicity we denote $\hat f = \ERM{\text{ERM}}^{0-1}$ and for all $f\in\cF$, $\cL_f^{0-1}=\cL_f$ and $R(f) = R^{0-1}(f)$.

It follows from the definition of the ERM that $P_N\cL_{\hat f}\leq 0$. Therefore, if we denote by $P_\cI$ (resp. $P_\cO$) the empirical measure supported on $\{(X_i, Y_i): i\in\cI\}$ (resp. $\{(X_i, Y_i): i\in\cO\}$), we have
\begin{align*}
 R(\hat f)-R(\bayes)&= (P-P_N)\cL_{\hat f}+P_N \cL_{\hat f}\leq (P-P_N)\cL_{\hat f} = \frac{|\cI|}{N}(P-P_\cI)\cL_{\hat f} + \frac{|\cO|}{N}(P-P_\cO)\cL_{\hat f}\\
& \le \frac{|\cI|}{N}\sup_{f\in \cF}(P-P_{\cI})\cL_f+\frac{2|\cO|}N
\end{align*}because $|\cL_{\hat f}|\leq 1$ a.s.. Then, by the bounded difference inequality \cite[Theorem 6.2]{concentration}, since all $f\in \cF$ satisfies $-1\le \cL_f\le 1$, one has, for any $x>0$,
\[
\bP\pa{\sup_{f\in \cF}(P-P_{\cI})\cL_f\ge \E[\sup_{f\in \cF}(P-P_{\cI})\cL_f]+x}\le e^{-2|\cI|x^2}\enspace.
\]


Furthermore, by the symmetrization argument (cf. Chapter~4 in \cite{MR2814399}),
\[
\E[\sup_{f\in \cF}(P-P_{\cI})\cL_f]\le 2\frac{\cR(\cL_\cF)}{|\cI|}\enspace.
\]
Therefore, for any $x>0$, with probability larger than $1-e^{-2|\cI|x^2}$,
\[
 R(\hat f)-R(\bayes)\le \frac{2\cR(\cL_\cF)}{N}+x+\frac{2|\cO|}{N}\enspace.
\]
The proof is completed by choosing $x=\sqrt{K/(2|\cI|)}$.

\subsection{Proof of Theorem~\ref{thm:MOM}}\label{sec:ProofThm1}
Let $\bayes\in \argmin_{f\in F}P\ell_f$. By definition, one has $\MOM{K}{\ell_{\ERM{\text{MOM},K}}}\le \MOM{K}{\ell_{f^*}} $, therefore,
\begin{equation}\label{eq:Vapnik}
R(\ERM{\text{MOM},K})-R(\bayes)\le P\ell_{\ERM{\text{MOM},K}}-\MOM{K}{\ell_{\ERM{\text{MOM},K}}} - \left(P\ell_{f^*}-\MOM{K}{\ell_{f^*}}\right)\enspace.
\end{equation} Let us now control the two expressions in the right-hand side of \eqref{eq:Vapnik}. Let $x>0$. We have
\begin{equation*}
\bP\left[P \ell_{f^*} - \MOM{K}{\ell_{f^*}}> x\right] = \bP \left[\sum_{k=1}^K I(P\ell_{f^*} - P_{B_k}\ell_{f^*}>x)\geq \frac{K}{2} \right] = \sum_{k=K/2}^K \binom{K}{k} p^k (1-p)^{K-k}\leq p^{K/2} 2^K
\end{equation*}where $p=\bP[P\ell_{f^*} - P_{B_k}\ell_{f^*}>x]$. Using Markov inequality together with ${\rm var}(\ell_{f^*})\leq 2 L^2 \E (f^*(X))^2\leq 2L^2 \theta^2$, we obtain
\begin{equation*}
\bP\left[P \ell_{f^*} - \MOM{K}{\ell_{f^*}}> x\right]\leq \left(\frac{4 {\rm var}(\ell_{f^*})K}{N x^2}\right)^{K/2}\leq \left(\frac{8L^2\theta^2K}{N x^2}\right)^{K/2}=\exp(-K/2)
\end{equation*}when $x=2L\theta \sqrt{2eK/N}$.

Now, for any $x>0$, one has $\sup_{f\in F}\MOM{K}{P\ell_f-\ell_f}>x$ iff 
\begin{equation}\label{eq:Obj}
\sup_{f\in F}\sum_{k=1}^KI\set{(P-P_{B_k})\ell_f>x}\ge \frac{K}2 \enspace. 
\end{equation}

Let us now control the probability that \eqref{eq:Obj} holds via an adaptation of the small ball method \cite{MR3431642,MR3367000}. Let $x>0$ and let $\phi(t)=(t-1)I\{1\le t\le 2\}+I\{t\ge 2\}$ be defined for all $t\in\R$. As $\phi(t)\ge I\{t\ge 2\}$, one has
\begin{multline*}
\sup_{f\in F}\sum_{k=1}^KI\set{(P-P_{B_k})\ell_f>x}\\
\le \sup_{f\in F}\sum_{k\in\cK}\E\cro{\phi\pa{2(P-P_{B_k})\ell_f/x}}+|\cO|+\sup_{f\in F}\sum_{k\in \cK}\pa{\phi\pa{2(P-P_{B_k})\ell_f/x}-\E\cro{\phi\pa{2(P-P_{B_k})\ell_f/x}}}  
\end{multline*}where we recall that $\cK=\{k\in\{1, \cdots, K\}: B_k\cap \cO = \emptyset\}\}$.

Since, $\phi(t)\le I\{t\ge 1\}$ and for all $f\in F$, $\text{Var}(\ell_f)\leq  2L^2\E f(X)^2\leq 2L^2\theta_2^2$, we have for all $f\in F$ and $k\in\cK$, 
\[
\E\cro{\phi\pa{2(P-P_{B_k})\ell_f/x}}\le \bP\pa{(P-P_{B_k})\ell_f\ge \frac x2}\le \frac{4\text{Var}(\ell_f)}{x^2|B_k|}\le \frac{8L^2\param{2}^2K}{x^2N}\enspace.
\]
One has therefore
\begin{multline*}
\sup_{f\in F}\sum_{k=1}^KI\set{(P-P_{B_k})\ell_f>x}\\
\le K\pa{\frac{8L^2\param{2}^2K}{x^2N}+\frac{|\cO|}K+\sup_{f\in F}\frac1K\sum_{k\in \cK}\pa{\phi\pa{\frac{2(P-P_{B_k})\ell_f}x}-\E\cro{\phi\pa{\frac{2(P-P_{B_k})\ell_f}x}}}}\enspace.  
\end{multline*}
As $0\le \phi(\cdot)\le 1$, by the bounded-difference inequality, for any $y>0$, with probability larger than $1-e^{-2y^2K}$,
\begin{align*}
\sup_{f\in F} \frac1K\sum_{k\in \cK}&\pa{\phi\pa{\frac{2(P-P_{B_k})\ell_f}x}-\E\cro{\phi\pa{\frac{2(P-P_{B_k})\ell_f}x}}}\\
& \le\E\cro{ \sup_{f\in F}\frac1K\sum_{k\in \cK}\pa{\phi\pa{\frac{2(P-P_{B_k})\ell_f}x}-\E\cro{\phi\pa{\frac{2(P-P_{B_k})\ell_f}x}}}}+y\enspace.
\end{align*}
Now, by the symmetrization inequality,
\begin{multline*}
\E\cro{\sup_{f\in F} \frac1K\sum_{k\in \cK}\pa{\phi\pa{\frac{2(P-P_{B_k})\ell_f}x}-\E\cro{\phi\pa{\frac{2(P-P_{B_k})\ell_f}x}}}}\\
\le 2\E\cro{ \sup_{f\in F}\frac1K\sum_{k\in \cK}\epsilon_k\phi\pa{\frac{2(P-P_{B_k})\ell_f}x}}\enspace. 
\end{multline*}
Since $\phi$ is $1$-Lipschitz and $\phi(0)=0$, by the contraction principle (see \cite[Chapter~4]{MR2814399} or more precisely equation~(2.1) in \cite{saintflour}),
\begin{equation*}
\E\cro{ \sup_{f\in F}\frac1K\sum_{k\in \cK}\epsilon_k\phi\pa{\frac{(P-P_{B_k})\ell_f}x}}\le \E\cro{ \sup_{f\in F}\frac{1}{xK}\sum_{k\in \cK}\epsilon_k(P-P_{B_k})\ell_f} \enspace.
\end{equation*}
By the symmetrization principle, 
\[
\E\cro{ \sup_{f\in F}\frac{2}{xK}\sum_{k\in \cK}\epsilon_k(P-P_{B_k})\ell_f} \le \frac{2}{xN}\E\cro{ \sup_{f\in F}\sum_{i\in \cJ}\epsilon_i\ell_f(X_i,Y_i)}\enspace.
\]
Finally, since $\ell$ is $L$-Lipschitz, by the contraction principle (see equation~(2.1) in \cite{saintflour}), 
\[
\E\cro{ \sup_{f\in F}\sum_{i\in \cJ}\epsilon_i\ell_f(X_i,Y_i)}\le 2L\cR(F)\enspace.
\]

Thus, for any $y>0$, with probability larger than $1-\exp(-2y^2K)$,
\begin{align*}
\sup_{f\in F}\sum_{k=1}^KI\set{(P-P_{B_k})\ell_f>x}
\le K\pa{\frac{8L^2\param{2}^2K}{x^2N}+\frac{|\cO|}K+y+\frac{4L\cR(F)}{xN}}\enspace.  
\end{align*}
Let $\Delta=1/4-|\cO|/K$ and let $y=\Delta$ and $x =  8L\max\left(\param{2}\sqrt{K/N}, 4\cR(F)/ N\right)$ so
\begin{align*}
\bP\pa{\sup_{f\in F}\sum_{k=1}^KI\set{(P-P_{B_k})\ell_f>x}
< \frac K2}\ge 1-e^{-\Delta^2K/8}\enspace.  
\end{align*}
Going back to \eqref{eq:Obj}, this means that
\[
\bP\pa{\sup_{f\in F}\MOM{K}{\ell_f-P\ell_f}\le 4L\max\left(\param{2}\sqrt{\frac{K}{N}}, \frac{4\cR(F)}{N}}\right)\ge 1-\exp(-2\Delta^2K)\enspace.  
\]
Plugging this result in \eqref{eq:Vapnik} concludes the proof of the theorem.

\subsection{Proof of Theorem~\ref{thm:convergence}}

First, hypothesis \ref{hyp:local_conv} implies that, for any $t\geqslant 1$, 
there exists an open convex neighbor $B$ of $u_{t-1}$ such that for all $u\in B$, $P_{B_{k_{med}}(\sigma_{t-1})}(\ell_{f_{u}})\in \MOM{K}{\ell_{f_{u}}}$. In particular therefore, for any $u\in B$,
\begin{equation}\label{eq:echange_deriv_MOM}
 \frac{K}{N}\sum_{b \in B_{k_{med}}(\sigma_{t-1})}\nabla_{u}\ell_{f_u}(X_b,Y_b)=\nabla_u\MOM{K}{f_{u}}.
 \end{equation}
From here, the rest of the proof is identical to the proof of consistency of the SGD. As a proof of concept, we reproduce here the arguments of \cite{SGD}. In what follows, all results are given conditionally on $\cD_N$. We set
\begin{equation*}
G_{t-1} = \nabla_u\MOM{K}{f_{u}}_{|u=u_{t-1}} = \frac{K}{N}\sum_{i \in B_{k_{med}}(\sigma_{t-1})}\nabla_u\ell_{f_{u}}(X_i,Y_i)_{|u=u_{t-1}}\enspace.
\end{equation*}

By definition of the update of the algorithm, we have 
\begin{align*}
\norm{u_t-u_{min}}_2^2 =& \norm{u_{t-1}-u_{min} -\eta_{t-1} G_{t-1}}_2^2\\
=&\norm{u_{t-1}-u_{min}}_2^2 -2\eta_{t-1} \inr{u_{t-1}-u_{min}, G_{t-1}} +\eta_t^2 \norm{G_{t-1}}_2^2.
\end{align*}
Let $F_t$ be the $\sigma$-algebra spanned by $\sigma_1,\dots,\sigma_{t-2}$, we have
\begin{align*}
\E\left[\norm{u_t-u_{min}}_2^2| F_t\right] =  \norm{u_{t-1}-u_{min}}_2^2 -2\eta_{t-1} \inr{ u_{t-1}-u_{min}, \E \left[G_{t-1}| F_t\right]} + \eta_{t-1}^2 \E \left[\norm{G_{t-1}}_2^2 | F_t\right].
\end{align*}
Then, we use Jensen's inequality, we have
\begin{multline}\label{eq:ineq_rec}
\E\left[\norm{u_t-u_{min}}_2^2| F_t\right] \le \norm{u_{t-1}-u_{min}}_2^2 -2\eta_{t-1} \frac{K}{N}\sum_{i \in B_{k_{med}}(\sigma_{t-1})}\inr{u_{t-1}-u_{min},\E\left[\left.\nabla_{u_{t-1}}\ell_{f_{u_{t-1}}}(X_i,Y_i)\right| F_t\right]} \\
+\eta_{t-1}^2\frac{K}{N}\sum_{i \in B_{k_{med}}(\sigma_{t-1})}\E\left[\left.\norm{\nabla_{u_{t-1}}\ell_{f_{u_{t-1}}}(X_i,Y_i) }_2^2\right|F_t\right].
\end{multline}
By convexity of $u\mapsto \ell_{f_u}(x,y)$, for any $x,y$, we have 
\[
\inr{u_{t-1}-u_{min},\nabla_{u_{t-1}}\ell_{f_{u_{t-1}}}(X_i,Y_i)}\geqslant \ell_{f_{u_{t-1}}}(X_i,Y_i)-\ell_{f_{u_{\min}}}(X_i,Y_i)
\]
Therefore,
\[
\frac{K}{N}\sum_{i \in B_{k_{med}}(\sigma_{t-1})}\inr{u_{t-1}-u_{min},\E\left[\left.\nabla_{u_{t-1}}\ell_{f_{u_{t-1}}}(X_i,Y_i)\right| F_t\right]}\geqslant \E\left[\MOM{K}{\ell_{f_{u_{t-1}}}}-P_{B_{k_{med}}(\sigma_{t-1})}\ell_{f_{u_{\min}}}|F_t\right]
\]
Hence, using hypothesis \ref{hyp:lipshitz} and  \ref{hyp:convex}, we get
\begin{equation}\label{eq:consistence_1}
\E\left[\norm{u_t-u_{min}}_2^2| F_t\right] \le  \norm{u_{t-1}-u_{min}}_2^2  +\eta_{t-1}^2L^2.
\end{equation}
Now, defining $h(t)=\norm{u_t-u_{min}}_2^2$ and $\delta_t=\1\{\E[h(t)-h(t-1)|F_t]>0\}$, it follows from \eqref{eq:consistence_1} that
$$\E[\delta(t)(h(t)-h(t-1))]=\E[\delta(t)\E[h(t)-h(t-1)|F_t]]\le \eta_t^2L^2 \enspace.$$
Using hypothesis \ref{hyp:seriesconv}, we have that
$$\sum_{t\geqslant 1}\E[\delta(t)(h(t)-h(t-1))]< \infty.$$
From Doob's first martingale convergence theorem (see for example \cite[Theorem 5.2.8]{martingale}), there exists $h_\infty \in \R_+$ such that 
$$h(t)\xrightarrow[t\to \infty]{a.s}h_\infty.$$
We plug this convergence in the inequality \eqref{eq:ineq_rec} to get 
\begin{align*}
2\eta_t (u_{t-1}-u_{min})^T\frac{K}{n}\E\left[\left.\sum_{i \in B_{k_{med}}(\sigma_{t-1})}\nabla_{u_{t-1}}\ell_{f_{u_{t-1}}}(X_i,Y_i)\right| F_t\right] &\le h(t-1)-\E[h(t)|F_t]+\eta_t^2L^2 \\
&\le h(t-1)+\eta_t^2L^2.
\end{align*}
Then, using hypothesis \ref{hyp:seriesconv}, we have
$$\sum_{t=1}^{\infty}\eta_t (u_{t-1}-u_{min})^T\frac{K}{N}\E\left[\left.\sum_{i \in B_{k_{med}}(\sigma_{t-1})}\nabla_{u_{t-1}}\ell_{f_{u_{t-1}}}(X_i,Y_i)\right| F_t\right] <\infty.$$
From the hypotheses \ref{hyp:seriesdiv} and \ref{hyp:convex}, we can then conclude that 
$$\norm{u_{t-1}-u_{min}}\xrightarrow[t\to \infty]{a.s}0. $$

\subsection{Proof of Proposition~\ref{prop:example_gradient}}
We denote by $B(1)(f),\dots,B(K)(f) $ the blocks such that the corresponding block empirical mean $P_{B(i)(f)}(f(X_1^n))$ are sorted: $P_{B(1)(f)}(f(X_1^n))\ge \dots\ge P_{B(K)(f)}(f(X_1^n))$. To simplify the notations we denote by $K/2=\lceil K/2 \rceil$ that correspond to the rank of the median when $K$ is odd, we suppose $K$ odd.\\

Let us take $\epsilon \in \R^N$, we show that for all $\varepsilon$ with $||\varepsilon||$ sufficiently small, we have $B_{(K/2)}(p)=B_{(K/2)}(p+t\varepsilon)$ for all $t\in [0,1]$.

To show that, it is sufficient to check that we have 
$$\forall 1\le k\le K-1, \forall t  \in [0,1],\quad P_{B(k)(p)}(f_{p+t\varepsilon})-P_{B(k+1)(p)}(f_{p+t\varepsilon})\ge0.$$
We decompose this difference in three parts,
\begin{align*}
P_{B(k)(p)}(f_{p+t\varepsilon})-P_{B(k+1)(p)}(f_{p+t\varepsilon})\ge&P_{B(k)(p)}(f_p)-P_{B(k+1)(p)}(f_p)\\
&-\left|P_{B(k)(p)}(f_p)-P_{B(k)(p)}(f(_{p+t\varepsilon})\right|\\
&-\left|P_{B(k+1)(p)}(f_{p+t\varepsilon})-P_{B(k+1)(p)}(f_p)\right|
\end{align*}
The two last terms are controlled by the Lipshitz property of $p\mapsto f_p$,
$$\forall t\in [0,1],\quad  P_{B(k)(p)}(f_{p+t\varepsilon})-P_{B(k+1)(p)}(f_{p+t\varepsilon})\ge P_{B(k)(p)}(f_p)-P_{B(k+1)(p)}(f_p)-2tL||\varepsilon||_2.$$
Then, we use order statistics and interquantile range to control the first term. 

To do that, we use Renyi representation to say that the order statistics of the blocks can be expressed as functions of the order statistics of exponential variables. Let $E_1,\dots,E_K$ be K i.i.d sample with a standard Exponential law, then using the function $U(t)= F_p^\leftarrow(1-1/t)$ where $F_p$ is the c.d.f of the random variable $P_{B_k}(f_p)$ to make the link between the sample $P_{B(1)(p)},\dots,P_{B(K)(p)}$ and the sample $E_1,\dots,E_K$, we have
$$\forall 1\le k\le K, \quad P_{B(k)(p)}(f_p) \sim U\left(\exp\left(\sum_{i=k}^K\frac{E_{(i)}}{i}\right)\right),$$

Then, because the hazard function is the inverse of the derivative of $U$ and is bounded by $1/m_{N,K}$, we have by the mean value inequality,
$$U\left(e^{\sum_{i=k}^K\frac{E_{(i)}}{i}}\right)-U\left(e^{\sum_{i=k+1}^{K}\frac{E_{(i)}}{i}}\right) \ge \frac{m_{N,K}E}{k}$$
hence, for all $\varepsilon' >0$ 
\begin{align*}
\P\left(P_{B(k)(p)}(f_p)-P_{B(k+1)(p)}(f_p) \le \varepsilon'\right)&\le \P\left( \frac{m_{N,K}E}{k}\le \varepsilon'\right)\\
&= 1-e^{-\frac{k\varepsilon'}{m_{N,K}}}.
\end{align*}
which we can rewrite as: for all $\delta \in (0,1)$, with probability greater than $\delta$, 
$$P_{B(k)(p)}(f_p)-P_{B(k+1)(p)}(f_p) \ge \frac{m\ln(1/\delta)}{k}.$$
Finally, with probability greater than $\delta$,
$$\forall t \in [0,1], \quad P_{B(k)(p)}(f_{p+t\varepsilon})-P_{B(k+1)(p)}(f_{p+t\varepsilon})\ge  \frac{m_{N,K}\ln(1/\delta)}{k}-2Lt||\varepsilon||_2.$$
Which implies that with probability greater than $\exp(-ktL||\varepsilon||_2/(2m_{N,K}))$,  
$$\forall t \in [0,1], \quad P_{B(k)(p)}(f_{p+t\varepsilon})-P_{B(k+1)(p)}(f_{p+t\varepsilon}) \ge 0.$$
Then, taking the union bound for $1\le k\le K-1$,  we have
$$\P\left(\forall t \in [0,1],\quad B_{(K/2)}(p)=B_{(K/2)}(p+t\varepsilon)\right)\ge 1-\sum_{k=1}^{K-1} 1-\exp(-ktL||\varepsilon||_2/(2m_{N,K})) $$

We want to compute the partial derivatives of $\psi_{f(p)}$, to do that we denote by $e_1,\dots,e_n \in \R^d$ the canonical basis of $\R^d$ and we define $\varepsilon_{n}^j=\frac{e_j}{2^n}m_{N,K}$. We denote by $A_n$ the event 
$$A_n^j=\left\{ \forall t \in [0,1],\quad , \forall n'\ge n, \quad B_{(K/2)}(p)=B_{(K/2)}(p+t\varepsilon_{n'}^j) \right\}$$ and we study the limiting event $\Omega^j=\overline{\lim}_{n\to\infty}A_n^j $.

First, let us note that for all $1\le j\le d$, the sequence of set $(A_n^j)_n$ is non-increasing, hence 
$$\Omega^j=\overline{\lim}_{n\to\infty}A_n^j=\underline{\lim}_{n\to\infty}A_n^j=(\overline{\lim}_{n\to\infty}(A_n^j)^c)^c,$$
then, for all $1\le j\le d$, we can study the $\overline{\lim}_{n\to\infty}(A_n^j)^c$ with Borel-Cantelli Lemma. Indeed, we have 
\begin{align*}
\P((A_n^j)^c)&\le\sum_{k=1}^K 1-\exp(-ktL||\varepsilon_n^j||_2/(2m_{N,K})) \\
&\le \sum_{k=1}^{K-1}\frac{ktL||\varepsilon_n^j||_2}{2m_{N,K}}\\
&= \frac{1}{2^n}\frac{K(K-1)tL}{4}
\end{align*}
Hence, the series $\sum_n \P((A_n^j)^c)$ converges and by Borel Cantelli Lemma, $\P(\overline{\lim}_{n\to\infty}(A_n^j)^c)=0$, then for all $1\le i\le d$, $\P(\Omega^j)=1$.

Rewriting that we have that 
$$\forall \omega \in \Omega^j,\quad \exists N \ge 1 : \quad \forall n \ge N, \quad \omega \in A_n^j .$$
Hence, $$\quad \exists N \ge 1 : \quad \forall n \ge N, \quad \forall t \in [0,1],\quad B_{(K/2)}(p)=B_{(K/2)}(p+t\varepsilon_n^j) ,$$
Hence, $$\quad \exists N \ge 1 : \quad \forall n \ge N, \quad \forall t \in [0,1],\quad B_{(K/2)}(p)=B_{(K/2)}(p+t\varepsilon_n^j) ,$$
which implies that for all $1\le j\le d$, $\forall n \ge N$,
\begin{align*}
\partial_j\psi_{f_p}(X)&=\lim_{t\to 0}\frac{\psi_{f_{p+t\epsilon_n^j}}(X)-\psi_{f_{p}}(X)}{t}\\
&=\lim_{t\to 0}\frac{P_{B(K/2)(p)}(f_{p+t\varepsilon_n^j})-P_{B(K/2)(p)}(f_{p})}{t}\\
&=\frac{1}{n/K}\lim_{t\to0}\sum_{i \in B(K/2)(p)} \frac{f_{p+t\varepsilon_n^j}(X_i)-f_{p}(X_i)}{t}\\
&=\frac{1}{n/K}\sum_{i \in B(K/2)(p)} \partial_jf_{p}(X_i)
\end{align*}

\section{Annex}

\subsection{Choice of the number of blocks}

\begin{figure}[h]
  \begin{center}

      \includegraphics[scale=0.7]{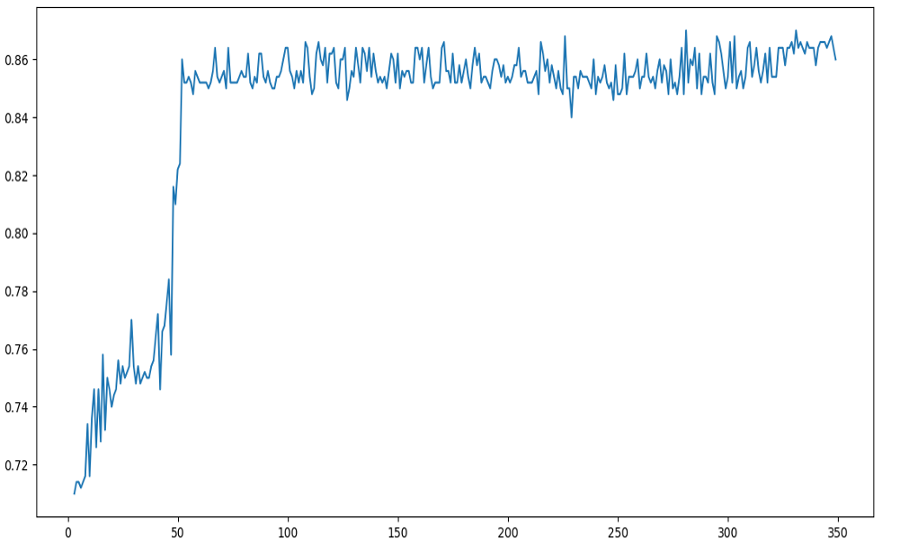}
 
    \caption{Plot of the accuracy on the toy dataset of Logistic Regression MOM as a function of $K$. \label{fig:plot_k}}
  \end{center}
\end{figure}

Let us study the behaviour of our algorithms when the number of blocks changes. We plot the accuracy as a function of $K$ averaged on 50 runs to have a good idea of the evolution of the performance with respect to $K$, the result is represented in figure \ref{fig:plot_k}.

\begin{figure}[h]
  \begin{center}

      \includegraphics[scale=0.7]{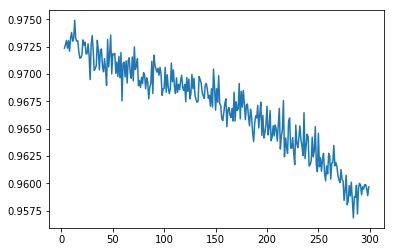}
 
    \caption{Plot of the accuracy on HTRU2 dataset of Logistic Regression MOM  as a function of $K$. \label{fig:plot_k_pulsar}}
  \end{center}
\end{figure}

There is a clear separation around $2|\cO|=60$ that is consistent with the theory. On the other hand the accuracy doesn't decrease when $K$ gets bigger one would expect. This may be due to the symmetry of the dataset.
If we run the same experiment on the real dataset, we get a much more regular plot, see Figure~\ref{fig:plot_k_pulsar}.

Figure \ref{fig:plot_k_pulsar}  confirms our predictions on clean datasets, the accuracy getting better as $K$ gets smaller (the MOM minimizer is the ERM when $K=1$ and ERM is optimal in the i.i.d. setup, \cite{LM13}). This may be due to the small number of outliers in this dataset.

\subsection{Illustration of convergence rate}
\begin{figure}
\begin{center}

    \subfloat[Scatter plot of the Moons dataset]{
      \includegraphics[width=0.4\textwidth]{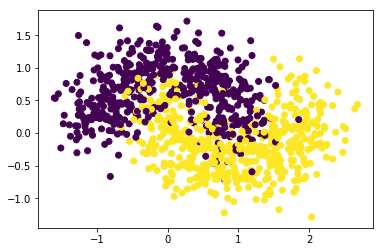}
      \label{sub:Moons}
                         }
    \subfloat[Scatter plot of the Gaussians dataset]{
      \includegraphics[width=0.4\textwidth]{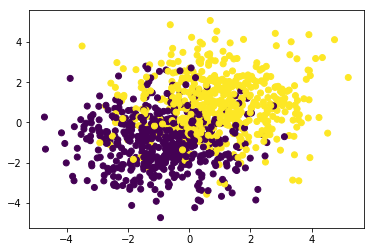}
      \label{sub:Gaussians}
                         }
    \caption{Scatter plot of the two dataset used in this section, the color represent the class of the points.}
    \label{fig:databases_rate}
\end{center}
\end{figure}

In this section, we estimate the rate of convergence of the MOM risk minimization algorithm Logistic Regression on two databases (see figure \ref{fig:databases_rate}). The first dataset is composed of points located on two interlaced half-circle with a Gaussian noise of standard deviation $0.3$, the two "moons" are each of a different class. We assume that these moons don't satisfy the margin property (we checked that the rate was slow for ERM algorithms, using the vanilla logistic regression). The second dataset is composed of two Gaussians $\mathcal{N}((-1,-1),1.4^2I_2) $ and $\mathcal{N}((1,1),1.4^2I_2)$ with respective label $1$ and $0$, we can prove that this dataset verifies the margin property needed to obtain fast rate in ERM

\begin{figure}
\begin{center}

    \subfloat[Convergence rate for the Moons dataset]{
      \includegraphics[width=0.45\textwidth]{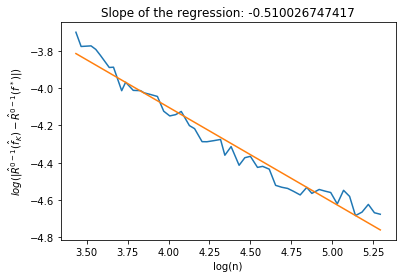}
      \label{sub:Moons}
                         }
    \subfloat[Convergence rate for the Gaussians dataset]{
      \includegraphics[width=0.45\textwidth]{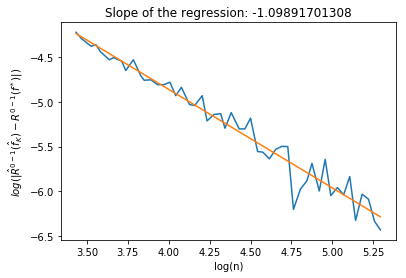}
      \label{sub:Gaussians}
                         }
    \caption{Plot of the logarithm of the excess risk as a function of $\log(n)$ in two cases: (a) where the margin assumption does not hold and (b) where the margin assumption holds. A linear regression is fitted on the curve, its slope is printed at the top of each figure revealing a slow $n^{-0.51}$ rate of convergence in case (a) and a fast $n^{-1.1}$ in case (b).}
    \label{fig:rates}
\end{center}
\end{figure}

There are no outliers in the datasets because we only want to test the rate of convergence.
To illustrate the rates of convergence of our algorithms, we plot the curve $\log\left(\left|\hat{R}^{0-1}(\hat{f} _K))-\hat{R}^{0-1}(f^*)\right|\right)$ as a function of $\log(n)$ where the risk is estimated by Monte-Carlo. 
The figure obtained for Logistic Regression MOM is represented in figure \ref{fig:rates}. It seems that MOM minimizers can achieve fast rates of convergence even if we did not prove them.

\begin{Remark}
We used random blocks sampled at each iteration for this application because it is the algorithm that we described earlier but even if we use one partition of blocks for the whole algorithm (as in the theory we developed) we obtain nonetheless fast rate for the Gaussians dataset.
\end{Remark}
\subsection{Comparison with robust algorithms based on M-estimators.}
\label{robust_classif_comparison}
In this section we compare the algorithm Logistic Regression MOM with two other algorithms based on M-estimators, these algorithms are studied on the toy dataset presented in Section \ref{simu}. 

One algorithm is a gradient on the Huber estimation of the loss function, it follows the same reasoning as MOM risk minimization and minimizes $\E[\l_f(X,Y)]$ using as a proxy the Huber estimator for this quantity. The Huber estimator is then defined as a $M$-estimator, denoted here $\hat\mu_f$, solution of 
$$\sum_{i=1}^n \psi_c(\hat\mu_f-\l(f(X_i),Y_i))=0 $$
where $\psi_c=\max(-c,\min(c,x))$ is the Huber function, $c>0$. Using this definition of $\hat \mu_f$, it is then easy to compute the gradient $\nabla \hat \mu_f$ and then use a gradient descent algorithm. The theory behind this algorithm is studied further in \cite{brownlees2015}.
 
The second algorithm uses a "redescending" loss function, in short we do ERM with a bounded loss function. Here we use Tukey biweight loss function rescaled by MADN scale estimator and IRLS algorithm to optimize the empirical risk.
\begin{figure}
\begin{center}
\includegraphics[scale=0.6]{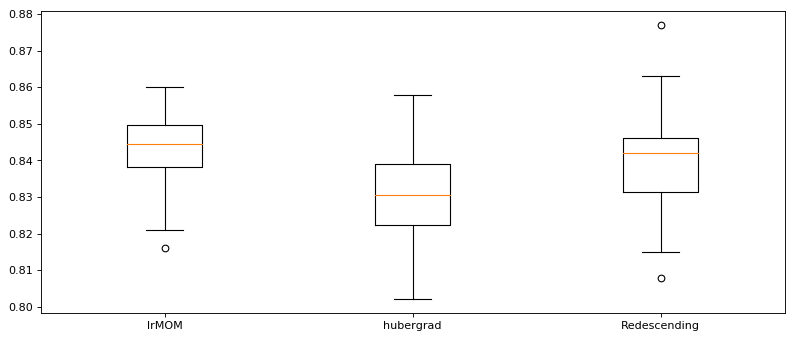}
\caption{Boxplot of the accuracy obtained on 50 training/test run (1000 training sample, $2\%$ corruption) of each algorithms on a 2-dimensional toy dataset.}
    \label{fig:comparison}

\end{center}
\end{figure}

\begin{figure}
\begin{center}

      \includegraphics[scale=0.6]{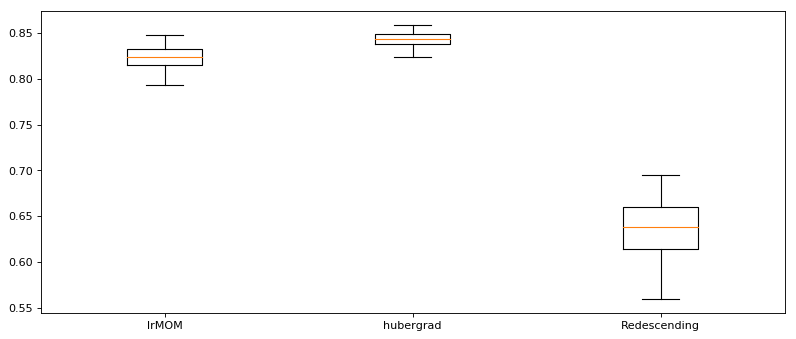}

    \caption{Boxplot of the accuracy obtained on 50 training/test run (2000 training sample, $2\%$ corruption) of each algorithms on a 200-dimensional toy dataset.}
    \label{fig:comparison_high_dim}
\end{center}
\end{figure}

Figure~\ref{fig:comparison} shows that all algorithms perform similarly on this easy, low dimensional dataset. The situation is quite different in higher dimension. In Figure \ref{fig:comparison_high_dim} we used a 200 dimensional dataset and the algorithm using a redescending loss function does not perform well. This may be due to local minima in which the algorithm gets stuck, as local minima are multiplied when the dimension gets higher. The other algorithms don't suffer this drawback since they use a ``projection by the loss function" that makes the problem one dimensional.

The algorithm using redescending loss functions is a simple gradient descent that has linear complexity. The Huber gradient algorithm estimates at each iteration a Huber estimator of location. The complexity of this estimator depends on the algorithm used but for most M-estimators a commonly used algorithm is an iteratively reweighted algorithm whose complexity is linear in the sample size. In practice we can nonetheless notice a great complexity of the Huber estimator in some cases where data are not well spread. In most cases, Logistic Regression MOM is the fastest among these three algorithms and the gradient Huber is the slowest, even though logistic regression may need a lot more iterations than the other algorithms. 

\begin{footnotesize}
\bibliographystyle{plain}
\bibliography{bibli}
\end{footnotesize}

\end{document}